\documentclass[11pt,reqno]{elsarticle}

\makeatletter
\def\ps@pprintTitle{%
 \let\@oddhead\@empty
 \let\@evenhead\@empty
 \def\@oddfoot{}%
 \let\@evenfoot\@oddfoot}
\makeatother

\usepackage{amsmath,amscd,amsfonts,amssymb,mathtools}
\usepackage{mathabx,breqn,multirow,comment,url}
\usepackage{booktabs,cool}
\usepackage{parskip}
\usepackage[mathscr]{eucal}
\usepackage[top=1.2in, bottom=1.2in, left=1in, right=1in]{geometry}
\usepackage{color}

\DeclareMathOperator{\cotan}{cot}

\DeclareMathOperator{\real}{Real}

\newtheorem{remark}{Remark}

\let\sep=\null

\providecommand{\e}[1]{\ensuremath{\times 10^{#1}}}

\makeatletter
\g@addto@macro\normalsize{%
  \setlength\abovedisplayskip{.4em}
  \setlength\belowdisplayskip{.4em}
  \setlength\abovedisplayshortskip{.4em}
  \setlength\belowdisplayshortskip{.4em}
}

\begin{document}

\begin{frontmatter}
\begin{abstract}
We describe a suite of fast algorithms for evaluating Jacobi polynomials, 
applying the corresponding discrete Sturm-Liouville eigentransforms
and calculating Gauss-Jacobi quadrature rules.  Our approach is based on
the well-known fact that 
Jacobi's differential equation admits a nonoscillatory phase function
which can be loosely approximated via an affine function over much of its
domain.
Our  algorithms perform  better than currently available methods in most respects.
We illustrate this with several numerical experiments, 
the source code for which is publicly available.

\end{abstract}


\title
{
Fast algorithms for Jacobi expansions via nonoscillatory phase functions
}

\author[jb]{James Bremer}
\ead{bremer@math.ucdavis.edu}

\address[jb]{Department of Mathematics, University of California, Davis}

\author[hz]{Haizhao Yang}
\ead{matyh@nus.edu.sg}
\address[hz]{Department of Mathematics, National University of Singapore}

\end{frontmatter}

\begin{section}{Introduction}

Expansions of functions in terms of orthogonal polynomials are widely
used in applied mathematics, physics, and numerical analysis.  And while highly effective
methods  for forming and manipulating expansions in terms of Chebyshev and Legendre polynomials are 
available,  there is substantial room for improvement in the analogous algorithms
for more general classes of Jacobi polynomials.

Here, we describe a collection of algorithms for forming and manipulating expansions
in Jacobi polynomials.   This includes methods for  
evaluating  Jacobi polynomials,  applying the forward and inverse Jacobi transforms,
and for calculating Gauss-Jacobi quadrature rules.
Our algorithms are based on the well-known observation that 
Jacobi's differential equation admits a nonoscillatory phase
function.
We combine this observation with a fast technique for computing nonoscillatory phase
functions and standard methods for exploiting the rank deficiency of matrices
which represent smooth functions to obtain our algorithms.

More specifically, given a pair of parameters $a$ and $b$,
both of which are the interval $\left(-\frac{1}{2},\frac{1}{2}\right)$, and a 
maximum degree of interest $N_{\mbox{\tiny max}} > 27$, we first numerically
construct nonoscillatory phase and amplitude functions $\psi^{(a,b)}(t,\nu)$
and $M^{(a,b)}(t,\nu)$ which represent the
Jacobi polynomials $P_\nu^{(a,b)}(x)$ of degrees between $27$ and $N_{\mbox{\tiny max}}$.
We evaluate the polynomials of lower degrees using the
well-known three-term recurrence relations; the lower bound of $27$ was chosen in light
of numerical experiments indicating it is close to optimal.
We restrict our attention to values of $a$ and $b$
in the interval $(-1/2,1/2)$ because,  in this regime, the solutions of 
Jacobi's differential equation are purely oscillatory and its phase
function is particularly simple.  When
$a$ and $b$ are larger, Jacobi's equation has turning points and the
corresponding phase function becomes more complicated.  The algorithms
of this paper continue to work for values of $a$ and $b$ modestly larger
than $1/2$, but they become less accurate as the parameters increase,
and eventually fail.    In a future work, we will discuss variants
of the methods described here which apply in the case of larger
values of the parameters (see the discussion in Section~\ref{section:conclusion}).

The relationship between $P_\nu^{(a,b)}$ and the nonoscillatory phase
and amplitude functions is
\begin{equation}
\tilde{P}_\nu^{(a,b)}(t) = M^{(a,b)}(t,\nu) \cos\left(\psi^{(a,b)}(t,\nu)\right),
\label{introduction:p}
\end{equation}
where $\tilde{P}_\nu^{(a,b)}$ is defined via the formula 
\begin{equation}
\tilde{P}_\nu^{(a,b)}(t) =
C_{\nu}^{a,b}\ P_\nu^{(a,b)}\left(\cos(t)\right) 
\sin\left(\frac{t}{2}\right)^{a+\frac{1}{2}}
\cos\left(\frac{t}{2}\right)^{b+\frac{1}{2}}
\label{introduction:ptilde}
\end{equation}
with
\begin{equation}
C_{\nu}^{a,b} = 
\sqrt{\left(2\nu + a + b + 1\right)\frac{\Gamma(1+\nu)\Gamma(1+\nu+a+b)}{\Gamma(1+\nu+a)
\Gamma(1+\nu+b)}}.
\label{introduction:C}
\end{equation}
The  constant (\ref{introduction:C}) ensures that the
$L^2\left(0,\pi\right)$ norm of $\tilde{P}_\nu^{(a,b)}$ is $1$
when $\nu$ is an integer; indeed,
the set $\left\{\tilde{P}_j^{(a,b)}\right\}_{j=0}^\infty$ is an orthonormal basis for $L^2\left(0,\pi\right)$.
The change of variables $x=\cos(t)$ is introduced 
because  it makes the singularities in the
phase and amplitude functions for Jacobi's differential equation more tractable.
We represent the functions $\psi^{(a,b)}$ and $M^{(a,b)}$
via their values at the nodes of a tensor product of piecewise Chebyshev grids.
Owing to the smoothness of the phase and amplitude functions, these representations
are quite efficient and can be constructed rapidly.  
Indeed, the asymptotic running time of our procedure for constructing
the nonoscillatory phase and amplitude functions
is  $\mathcal{O}\left(\log^2\left(N_{\mbox{\tiny max}}\right)\right)$.

Once the nonoscillatory phase and amplitude functions have been constructed,
the function $P_\nu^{(a,b)}$ can be evaluated at any point
$x \in (-1,1)$ and for any $27 \leq \nu \leq N_{\mbox{\tiny max}}$ in time  which is  independent
of $N_{\mbox{\tiny max}}$, $\nu$ and $x$ via (\ref{introduction:p}).
One downside of our approach is that the error which occurs when
Formula~(\ref{introduction:p}) is used to evaluate a Jacobi polynomial
numerically
 increases with the magnitude of the phase
function $\psi^{(a,b)}$ 
owing to the well-known difficulties in evaluating trigonometric
functions of large arguments. 
 This is not surprising, and the resulting
errors are in line with the 
condition number of evaluation of the function $P_\nu^{(a,b)}$.  
Moreover,  this phenomenon is hardly limited to the approach of this paper
and can be observed, for instance, when asymptotic expansions are used
to evaluate Jacobi polynomials (see, for instance, \cite{Bogaert-Michiels-Fostier} for an extensive
discussion of this issue in the case of asymptotic expansions
for Legendre polynomials).

Aside from applying the Jacobi transform and evaluating Jacobi polynomials,
nonoscillatory phase functions are also highly useful for calculating
the roots of special functions.
From (\ref{introduction:p}), we see that the roots of 
$P_\nu^{(a,b)}$ occur when
\begin{equation}
\psi^{(a,b)}(t,\nu) = \frac{\pi}{2} + n \pi
\end{equation}
with $n$ an integer.   Here,  we describe a method
for rapidly computing Gauss-Jacobi quadrature rules which exploits this fact.
The  $n$-point Gauss-Jacobi quadrature rule corresponding to the
parameters $a$ and $b$ is, of course, the unique quadrature rule of the form
\begin{equation}
\int_{-1}^1 f(x) (1-x)^a (1+x)^b\ dx 
\approx
\sum_{j=1}^n f(x_j) w_j,
\label{introduction:quadrule}
\end{equation}
which is exact when $f$ is a polynomial of degree less than or equal to $2n-1$
(see, for instance, Chapter~15 of \cite{Szego} 
for a detailed discussion of Gaussian quadrature rules).
Our algorithm can also  produce what we call the modified
Gauss-Jacobi quadrature rule.  That is, the $n$-point quadrature rule
of the form
\begin{equation}
\int_{0}^\pi f(\cos(t)) \cos^{2a+1}\left(\frac{t}{2}\right)   \sin^{2b+1}\left(\frac{t}{2}\right)\ dt 
\approx
\sum_{j=1}^n f(\cos(t_j)) 
\cos^{2a+1}\left(\frac{t}{2}\right)   \sin^{2b+1}\left(\frac{t}{2}\right)
w_j
\label{introduction:modrule}
\end{equation}
which is exact whenever $f$ is a polynomial of degree less than or equal to $2n-1$.  
The modified Gauss-Jacobi quadrature rule integrates products
of the functions $\tilde{P}_\nu^{(a,b)}$ of degrees between $0$ and $n-1$.
A variant of this algorithm which allows for the calculation of zeros of much more general classes
of special functions was previously published in \cite{BremerZeros}; however,
the version described here is specialized to the case of Gauss-Jacobi
quadrature rules and is considerably more efficient.  

We also describe a method for applying the Jacobi transform
using the nonoscillatory phase and amplitude functions.
For our purposes, the $n^{th}$ order discrete forward Jacobi transform
 consists of calculating the vector  of values
\begin{equation}
\left(
\begin{array}{c}
f(t_1) \sqrt{w_1} \\
f(t_2) \sqrt{w_2} \\
\vdots\\
f(t_n) \sqrt{w_n} \\
\end{array}
\right)
\label{introduction:transout}
\end{equation}
given the vector
\begin{equation}
\left(
\begin{array}{c}
\alpha_1 \\
\alpha_2\\
\vdots \\
\alpha_n
\end{array}
\right)
\label{introduction:transin}
\end{equation}
of the coefficients in the expansion
\begin{equation}
f(t) = \sum_{j=0}^{n-1} \alpha_j \tilde{P}_j^{(a,b)}(t);
\label{introduction:transexp}
\end{equation}
here,
 $t_1,\ldots,t_n,w_1\ldots,w_n$ are the nodes and weights
of the $n$-point modified Gauss-Jacobi quadrature rule corresponding
to the parameters $a$ and $b$.
The properties of this quadrature rule and the
weighting by square roots in (\ref{introduction:transout}) ensure that the
$n \times n$ matrix $\mathcal{J}_n^{(a,b)}$
which takes (\ref{introduction:transin}) to (\ref{introduction:transout}) is orthogonal.
It follows from (\ref{introduction:p}) 
 that the $(j,k)$ entry of $\mathcal{J}_n^{(a,b)}$ 
is
\begin{equation}
M^{(a,b)}(t_j,k-1) \cos\left( \psi^{(a,b)}(t_j,k-1) \right) \sqrt{w_j}.
\label{introduction:jn}
\end{equation}

Our method for applying $\mathcal{J}^{(a,b)}_n$, which is inspired
by the algorithm of \cite{Candes-Demanet-Ying} for applying Fourier integral transforms 
and a special case of that in \cite{Yang} for general oscillatory integral transforms,  exploits 
the fact that the matrix whose $(j,k)$  entry is 
\begin{equation}
 M^{(a,b)}(t_j, k-1)\exp\left( i \left(\psi^{(a,b)} (t_j, k-1) - (k-1) t_j\right) \right)\sqrt{w_j}
\label{introduction:matrix}
\end{equation}
is low rank. 
Since $\mathcal{J}_n^{(a,b)}$ is the real part of the Hadamard product
of the matrix defined by (\ref{introduction:matrix})
and the  $n \times n$ nonuniform FFT matrix whose $(j,k)$ entry is 
\begin{equation}
\exp\left( i (k-1) t_j\right),
\label{introduction:nufft}
\end{equation}
it can be applied through a small number of 
discrete nonuniform Fourier transforms.  This can be done quickly
via any number of fast algorithms for the nonuniform
Fourier transform 
(examples of such algorithms
include \cite{Greengard-NUFFT} and \cite{TownsendNUFFT}).
By ``Hadamard product,'' we mean the operation of  
entrywise multiplication of two matrices.    We will denote the Hadamard
product of the matrices $A$ and $B$ via $A \otimes B$ in this article.

We do not prove a rigorous bound on the rank of the matrix (\ref{introduction:matrix}) here; instead,
we offer the following heuristic argument.  For points $t$ bounded away from
the endpoints of the interval $(0,\pi)$, we have the following asymptotic approximation of the 
nonoscillatory phase function:
\begin{equation}
\psi^{(a,b)}(t,\nu) =  \left(\nu + \frac{a+b+1}{2}\right) t +c + 
\mathcal{O}\left(\frac{1}{\nu^2}\right)\ \ \mbox{as} \ \ \nu \to \infty
\label{introduction:asym}
\end{equation}
with $c$ a constant in the interval $(0,2\pi)$.  We prove this
in Section~\ref{section:jacobiphase}, below.
 From (\ref{introduction:asym}), we expect that the function
\begin{equation}
\psi^{(a,b)}(t,\nu) - t \nu
\label{introduction:psiapprox}
\end{equation}
is of relatively small magnitude, and hence that the rank of 
the matrix (\ref{introduction:matrix}) is low.
This is borne out by our numerical experiments
(see  Figure~\ref{figure:transform2} in particular).

We choose to approximate  $\psi^{(a,b)}(t,\nu)$ via
the linear function $\nu t$  rather than with a more complicated function
(e.g., a nonoscillatory combination of Bessel functions)
because this approximation leads to a method for
applying the Jacobi transform through repeated applications
of the fast Fourier transform, if the non-uniform Fourier 
transform in the Hadamard product just introduced is carried 
out via the algorithm in \cite{TownsendNUFFT}. It might be 
more efficient to apply the non-uniform FFT algorithm in 
\cite{Greengard-NUFFT}, depending on different computation 
platforms and compilers; but we do not aim at optimizing our 
algorithm in this direction. 
Hence, our algorithm  for applying the Jacobi transform
requires $r$ applications of the fast Fourier transform, 
where  $r$ is the
numerical rank of a matrix which is strongly related to (\ref{introduction:matrix}).
This makes  its  asymptotic complexity  $\mathcal{O}\left(r n \log(n)\right)$.
Again, we do not prove a rigorous bound on $r$ here.  However,
in the experiments of Section~\ref{section:experiments:transform}
we compare  our algorithm for applying the Jacobi transform
with Slevinsky's algorithm \cite{Slevinsky1} for applying
the Chebyshev-Jacobi transform.  The latter's running time is
$\mathcal{O}\left(n \log^2 (n) / \log(\log(n))\right)$, and the behavior
of our algorithm is similar.  This leads  us to  conjecture that 
\begin{equation}
r = \mathcal{O}\left( \frac{\log(n)}{\log(\log(n))}\right).
\end{equation}
Before our algorithm for applying the Jacobi
transform can be used, a precomputation
in addition to the construction
of the nonoscillatory phase and amplitude functions must be performed.  Fortunately,
the running time of this  precomputation procedure 
 is quite small.  Indeed,  for most values of $n$, it requires less time than
 a single application of the Jacobi-Chebyshev transform via \cite{Slevinsky1}.
Once the precomputation phase has been dispensed with, we our algorithm
for the application of the Jacobi transform is roughly ten times faster
than the algorithm of  \cite{Slevinsky1}.

The remainder of this article is structured as follows.  
Section~\ref{section:preliminaries}
summarizes a number of mathematical and numerical facts to be used in the rest of the paper.
In Section~\ref{section:precomp}, we detail the scheme we use to construct
the phase and amplitude functions.  Section~\ref{section:quad} describes
our method for the calculation of Gauss-Jacobi quadrature rules. 
In Section~\ref{section:transform}, we give our
algorithm for applying the forward and inverse Jacobi transforms. 
In Section~\ref{section:experiments}, we describe the results of numerical 
experiments conducted to assess the performance of the algorithms of this paper.
Finally, we conclude in Section~\ref{section:conclusion} with a few brief remarks
regarding this article and a discussion of directions for further work.

\end{section}

\begin{section}{Preliminaries}


\begin{subsection}{Jacobi's differential equation and the solution of the second kind}

The Jacobi function of the first kind is given by the formula
\begin{equation}
P_\nu^{(a,b)}(x) = \frac{\Gamma(\nu+a+1)}{\Gamma(\nu+1)\Gamma(a+1)} 
\Hypergeometric{2}{1}{-\nu,\nu+a+b+1}{a+1}{\frac{1-x}{2}},
\label{jacobieq:p}
\end{equation}
where $\Hypergeometric{2}{1}{a,b}{c}{z}$ denotes Gauss' hypergeometric function
(see, for instance, Chapter~2 of \cite{HTFI}).
We also define the Jacobi function of the second kind via 
\begin{dmath}
Q_\nu^{(a,b)} = \cotan(a \pi) P_\nu^{(a,b)}(x) 
-  \\
\frac{ 2^{a+b}\Gamma(\nu+b+1)\Gamma(a)}{\pi \Gamma(\nu+b+a+1)}
(1-x)^{-a} (1+x)^{-b} 
\Hypergeometric{2}{1}{\nu+1,-\nu-a-b}{1-a}{\frac{1-x}{2}}.
\label{jacobieq:q}
\end{dmath}
Formula (\ref{jacobieq:p}) is valid for all values
of the parameters $a$, $b$, and $\nu$ and all values
of the argument $x$ of interest to us here.  
Formula (\ref{jacobieq:q}), on the other hand, is invalid
when $a=0$ (and more generally when $a$ is an integer).  However, the limit as $a \to 0$ of $Q_\nu^{(a,b)}$ exists and various
series expansions for it can be derived rather easily
(using, for instance, the well-known results regarding hypergeometric
series which appear in Chapter~2 of \cite{HTFI}).
Accordingly, we view (\ref{jacobieq:q}) as defining
$Q_\nu^{(a,b)}$ for all $a,b$ in our region of interest $(-1/2,1/2)$.

Both $P_\nu^{(a,b)}$ and $Q_\nu^{(a,b)}$ are solutions of Jacobi's differential equation 
\begin{equation}
(1-x^2) \psi''(x)  + (b-a - (a+b+2)x)\psi'(x) + \nu(\nu+a+b+1)\psi(x) = 0
\label{jacobieq:1}
\end{equation}
(see, for instance, Section 10.8 of \cite{HTFII}).
We note that our choice of $Q_\nu^{a,b}$ is nonstandard
and differs from the convention used in \cite{Szego} and \cite{HTFII}.   However,
it is natural in light of Formulas~(\ref{asym:p}) and (\ref{asym:q}),
(\ref{hahn:p}) and (\ref{hahn:q}), and especially (\ref{jacobiphase:Mapprox}),
all of which appear below.

It can be easily verified that the functions
$\tilde{P}_\nu^{(a,b)}$ and $\tilde{Q}_\nu^{(a,b)}$ 
defined via (\ref{introduction:ptilde}) and the formula
\begin{equation}
\tilde{Q}_\nu^{(a,b)}(t) = 
C_{\nu}^{a,b} \ Q_\nu^{(a,b)}\left(\cos(t)\right) 
\sin\left(\frac{t}{2}\right)^{a+\frac{1}{2}}
\cos\left(\frac{t}{2}\right)^{b+\frac{1}{2}}
\label{preliminaries:qtilde}
\end{equation}
satisfy the second order differential equation
\begin{equation}
y''(t) + q^{(a,b)}_\nu(t) y(t) = 0,
\label{jacobieq:mod}
\end{equation}
where
\begin{dmath}
q_\nu^{(a,b)}(t) = \nu  (a+b+\nu +1)+\frac{1}{2} \csc ^2(t) (-(b-a) \cos (t)+a+b+1)-\frac{1}{4} ((a+b+1)
   \cot (t)-(b-a) \csc (t))^2.
\label{jacobieq:coef}
\end{dmath}
We refer to Equation~(\ref{jacobieq:mod}) as the modified Jacobi differential
equation.

\end{subsection}

\begin{subsection}{The asymptotic approximations  of Baratella and Gatteschi}

In \cite{Baratella-Gatteschi},   it is shown that there exist
sequences of real-valued functions 
$\left\{A_j\right\}$ and $\left\{B_j\right\}$ such
that for each nonnegative integer $m$,
\begin{equation}
\tilde{P}_\nu^{(a,b)}(t) =
C_\nu^{(a,b)} \frac{ p^{-a}}{\sqrt{2}} \frac{\Gamma(\nu+a+1)}{\Gamma(\nu+1)}
\left(
\sum_{j=0}^m \frac{A_j(t)}{p^{2j}} t^{\frac{1}{2}} J_a(pt) 
+
\sum_{j=0}^{m-1} \frac{B_j(t)}{p^{2j+1}} t^{\frac{3}{2}} J_{a+1}(pt) 
+ \epsilon_m(t,p)
\right),
\label{asym:p}
\end{equation}
where $J_\mu$ denotes the Bessel function of the first kind
of order $\mu$, $C_\nu^{(a,b)}$ is defined by  (\ref{introduction:C}), 
$p = \nu + (a+b+1)/2$, and 
\begin{equation}
\epsilon_m(t,p) = \begin{cases}
t^{a+\frac{5}{2}} \mathcal{O}\left(p^{-2m+a}\right), & 0 <  t \leq \frac{c}{p} \\
t \mathcal{O}\left(p^{-2m-\frac{3}{2}}\right), & \frac{c}{p} \leq t \leq \pi-\delta\\
\end{cases}
\label{asym:errorterm}
\end{equation}
with $c$ and $\delta$ fixed constants.
The first few coefficients in (\ref{asym:p}) are  given by 
$A_0(t) = 1$, 
\begin{equation}
B_0(t) = \frac{1}{4t} g(t),
\end{equation}
and
\begin{equation}
A_1(t) =  \frac{1}{8} g'(t) - \frac{1+2a}{8t} g(t) - \frac{1}{32} (g(t))^2 + 
\frac{a}{24} \left(3b^2 + a^2 - 1\right),
\end{equation}
where
\begin{equation}
g(t) = \left(\frac{1}{4}-a^2\right) 
\left( \cotan\left(\frac{t}{2}\right)-\frac{2}{t}\right)
- \left(\frac{1}{4}-b^2\right) \tan\left(\frac{t}{2}\right).
\end{equation}
Explicit formulas for the higher order coefficients are not known,
but formulas defining them are given in \cite{Baratella-Gatteschi}.
The analogous expansion of the function of the
second kind is
\begin{equation}
\tilde{Q}_\nu^{(a,b)}(t) \approx
C_\nu^{(a,b)} \frac{p^{-a}}{\sqrt{2}} \frac{\Gamma(\nu+a+1)}{\Gamma(\nu+1)}
\left(
\sum_{s=0}^m \frac{A_s(t)}{p^{2s}} t^{\frac{1}{2}} Y_a(pt) 
+
\sum_{s=0}^{m-1} \frac{B_s(t)}{p^{2s+1}} t^{\frac{3}{2}} Y_{a+1}(pt) 
+ \epsilon_m(t,p)\right).
\label{asym:q}
\end{equation}

An alternate asymptotic expansion of $P_\nu^{(a,b)}$ whose form is similar to
(\ref{asym:p}) appears in \cite{Frenzen-Wong},
and several more general Liouville-Green type expansions for Jacobi polynomials
can be found in   \cite{Dunster-Jacobi}.  The expansions of \cite{Dunster-Jacobi}
involve a fairly complicated change of variables which complicates  
their use in numerical algorithms.  

\label{section:gatteschi}
\end{subsection}

\begin{subsection}{Hahn's trigonometric expansions}

The asymptotic formula
\begin{equation}
\tilde{P}_\nu^{(a,b)}(t) = 
C_{\nu}^{(a,b)} \frac{2^{2p}}{\pi} B(\nu+a+1,\nu+b+1)
\sum_{m=0}^{M-1} \frac{ f_m(t) }{2^m \left(2p+1\right)_m} + R_{\nu,m}^{(a,b)}(t)
\end{equation}
appears in a slightly different form in \cite{Hahn}.   Here,
 $(x)_m$ is the Pochhammer symbol, $B$ is the beta function,
$p$ is once again equal to $\nu + (a+b+1)/2$,  and
\begin{equation}
f_m(t) = \sum_{l=0}^m 
\frac{\left(\frac{1}{2}+a\right)_l \left(\frac{1}{2}-a\right)_l
\left(\frac{1}{2}+b\right)_{m-l} \left(\frac{1}{2}-b\right)_{m-l}}
{\Gamma(l+1) \Gamma(m-l+1)}  
\frac
{\cos\left(\frac{1}{2}(2p + m)t - \frac{1}{2} \left(a+l+\frac{1}{2}\right)\pi\right)}
{\sin^l\left(\frac{t}{2}\right)\cos^{m-l}\left(\frac{t}{2}\right)}.
\label{hahn:p}
\end{equation}
The remainder term $R_{\nu,m}^{(a,b)}$ is bounded by twice the magnitude of the first
neglected term when $a$ and $b$ are in the interval $\left(-\frac{1}{2},\frac{1}{2}\right)$. 
The analogous approximation of the Jacobi function of the second kind
is
\begin{equation}
\tilde{Q}_\nu^{(a,b)}(t) \approx
C_{\nu}^{(a,b)} \frac{2^{2p}}{\pi} B(\nu+a+1,\nu+b+1)
\sum_{m=0}^{M-1} \frac{ g_m(t) }{2^m \left(2p+1\right)_m},
\end{equation}
where
\begin{equation}
g_m(t) = \sum_{l=0}^m 
\frac{\left(\frac{1}{2}+a\right)_l \left(\frac{1}{2}-a\right)_l
\left(\frac{1}{2}+b\right)_{m-l} \left(\frac{1}{2}-b\right)_{m-l}}
{\Gamma(l+1) \Gamma(m-l+1)}  
\frac
{\sin\left(\frac{1}{2}(2p + m)t - \frac{1}{2} \left(a+l+\frac{1}{2}\right)\pi\right)}
{\sin^l\left(\frac{t}{2}\right)\cos^{m-l}\left(\frac{t}{2}\right)}.
\label{hahn:q}
\end{equation}
\label{section:hahn}
\end{subsection}

\begin{subsection}{Some elementary facts regarding phase functions for second order differential equations}
We say that $\psi$ is a phase function for the second order differential
equation
\begin{equation}
y''(t) + q(t) y(t) = 0 \ \ \mbox{for all} \ \ \sigma_1 \leq t \leq \sigma_2
\label{phase:diffeq}
\end{equation}
provided the functions
\begin{equation}
 \frac{\cos\left(\psi(t)\right)}{\sqrt{\left|\psi'(t)\right|}}
\label{phase:u1}
\end{equation}
and
\begin{equation}
\frac{\sin\left(\psi(t)\right)}{\sqrt{\left|\psi'(t)\right|}}
\label{phase:v1}
\end{equation}
form a basis in its space of solutions.   By repeatedly
differentiating (\ref{phase:u1}) and (\ref{phase:v1}),
it can be shown that  $\psi$ is a phase
function for (\ref{phase:diffeq}) if and only if 
its derivative satisfies the nonlinear ordinary differential equation
\begin{equation}
 q(t) 
- (\psi'(t))^2
- \frac{1}{2}\left(\frac{\psi'''(t)}{\psi'(t)}\right)
+ \frac{3}{4}
\left(\frac{\psi''(t)}{\psi'(t)}\right)^2 = 0
\ \ \mbox{for all}\  \sigma_1 \leq t \leq \sigma_2.
\label{phase:kummer}
\end{equation}

A pair $u$, $v$ of real-valued,  linearly independent solutions
of (\ref{phase:diffeq})  determine a phase function $\psi$ up to an integer
multiple of $2\pi$.  Indeed, it can be easily seen that the function
\begin{equation}
\psi'(t) = \frac{W}{(u(t))^2+(v(t))^2},
\label{phase:psip}
\end{equation}
where $W$ is the necessarily constant Wronskian  of the pair $\{u,v\}$, satisfies
(\ref{phase:kummer}).  As a consequence, if we define
\begin{equation}
\psi(t) = C + \int_{\sigma_1}^t 
 \frac{W}{(u(s))^2+(v(s))^2} \ ds
\end{equation}
with $C$ an appropriately chosen constant, then
\begin{equation}
u(t) = \sqrt{W}\ \frac{\cos\left(\psi(t)\right)}{\sqrt{\left|\psi'(t)\right|}}
\label{phase:u2}
\end{equation}
and
\begin{equation}
v(t) = \sqrt{W}\ \frac{\sin\left(\psi(t)\right)}{\sqrt{\left|\psi'(t)\right|}}.
\label{phase:v2}
\end{equation}
The requirement that (\ref{phase:u2}) and (\ref{phase:v2}) hold clearly 
determines  the value of $\mbox{mod}(C,2\pi)$.  

If $\psi$ is a phase function for (\ref{phase:diffeq}),
then we refer to the function
\begin{equation}
M(t) = \sqrt{\frac{W}{\left|\psi'(t)\right|}}
\label{phase:mpsip}
\end{equation}
as the corresponding amplitude function.  
%
Though a straightforward computation
it can be verified that the 
the square  $N(t) = \left(M(t)\right)^2$ of the amplitude function
satisfies the third order linear ordinary differential equation
\begin{equation}
N'''(t) + 4 q(t) N'(t) + 2 q'(t) N(t) = 0
\ \ \mbox{for all}\ \ \sigma_1 < t < \sigma_2.
\label{phase:Neq}
\end{equation}

\label{section:phase}
\end{subsection}


\begin{subsection}{The nonoscillatory phase function for Jacobi's differential equation}

We will denote by $\psi^{(a,b)}(t,\nu)$ a phase function for 
the modified Jacobi differential equation (\ref{jacobieq:mod}) which, for each $\nu$,
gives rise to the  pair of solutions $\tilde{P}_\nu^{(a,b)}$ and $\tilde{Q}_\nu^{(a,b)}$.
We let $M^{(a,b)}(t,\nu)$ denote the corresponding amplitude function, so that 
\begin{equation}
\tilde{P}_\nu^{(a,b)}(t) = M^{(a,b)}(t,\nu) \cos\left(\psi^{(a,b)}(t,\nu) \right)
\label{jacobiphase:p}
\end{equation}
and
\begin{equation}
\tilde{Q}_\nu^{(a,b)}(t) = M^{(a,b)}(t,\nu) \sin\left(\psi^{(a,b)}(t,\nu) \right).
\label{jacobiphase:q}
\end{equation}
We use the notations $\psi^{(a,b)}(t,\nu)$ and  $M^{(a,b)}(t,\nu)$ rather than
$\psi_\nu^{(a,b)}(t)$ and $M_\nu^{(a,b)}(t)$, which might seem more natural, to emphasize that
the representations of the phase and amplitude functions we construct in this
article allow for their evaluation for a range of values of $\nu$ and $t$, but
only for particular fixed values of $a$ and $b$.

Obviously,
\begin{equation}
\left(M^{(a,b)}(t,\nu)\right)^2 = 
\left(P_\nu^{(a,b)}(t)\right)^2+\left(Q_\nu^{(a,b)}(t)\right)^2.
\label{jacobiphase:M}
\end{equation}
By replacing the Jacobi functions in (\ref{jacobiphase:M})
with the first terms in the expansions
Formulas~(\ref{asym:p})  and (\ref{asym:q}), we see that
for all $t$ in an interval of the form $(\delta,\pi-\delta)$ 
with $\delta$ a small positive constant and
all $a,b  \in (-1/2,1/2)$,
\begin{equation}
\begin{aligned}
\left(M^{(a,b)}(t,\nu)\right)^2
\sim
\left(C_\nu^{(a,b)} \frac{ p^{-a}}{\sqrt{2}} \frac{\Gamma(\nu+a+1)}{\Gamma(\nu+1)}\right)^2
t \left(
\left(J_{a}(pt)\right)^2 +  \left(Y_{a}(pt)\right)^2\right)
+ \mathcal{O}\left(\frac{1}{\nu^2}\right)
\ \ \mbox{as} \ \ \nu \to \infty
\end{aligned}
\label{jacobiphase:Mapprox}
\end{equation}
It is well known that the function
\begin{equation}
\left(J_{a}(t)\right)^2 +  \left(Y_{a}(t)\right)^2
\label{jacobiphase:bessel}
\end{equation}
is  nonoscillatory.    Indeed, it is completely
monotone on the interval $(0,\infty)$, as can be shown using
Nicholson's formula
\begin{equation}
\left(J_a(t)\right)^2
+
\left(Y_a(t)\right)^2
 = \frac{8}{\pi^2} \int_0^\infty K_0(2t\sinh(s)) \cosh(2a s)\ ds
\ \ \mbox{for all}\ \  t >0.
\label{jacobiphase:nicholson}
\end{equation}
A derivation of (\ref{jacobiphase:nicholson})  can be found in Section~13.73 of \cite{Watson};
that (\ref{jacobiphase:bessel}) is completely monotone follows easily from (\ref{jacobiphase:nicholson})
(see, for instance, \cite{Merkle}).
When higher order approximations of $\tilde{P}_{\nu}^{(a,b)}$ and $\tilde{Q}_{\nu}^{(a,b)}$
are inserted into (\ref{jacobiphase:M}), the resulting terms are all
 nonoscillatory combinations of Bessel functions.
From this it is clear that $M_{\nu}^{(a,b)}$ is nonoscillatory,
and Formula~(\ref{jacobiphase:Mapprox}) then implies that
 $\psi_\nu^{(a,b)}$ is also a nonoscillatory function.

It is well-known that
\begin{equation}
\left(J_a(t)\right)^2
+
\left(Y_a(t)\right)^2
\sim
\frac{2}{\pi}
\left(
\frac{1}{t }+ 
\frac{4a^2-1}{8t^2} + \frac{3(9-40a^2+16a^4)}{128t^5} + \cdots \right)
\ \ \mbox{as} \ \ t \to \infty;
\label{jacobiphase:bessel2}
\end{equation}
see, for instance, Section~13.75 of \cite{Watson}.  
From this, (\ref{jacobiphase:Mapprox}) and 
(\ref{introduction:C}),
we easily conclude that for $t$ in $(\delta,\pi-\delta)$
and $(a,b) \in (-1/2,1/2)$,
\begin{dmath}
\left(M^{(a,b)}(t,\nu)\right)^2 \sim
\frac{2p^{-2a}}{\pi}
\frac{\Gamma(1+\nu+a)\Gamma(1+\nu+a+b)}{\Gamma(1+\nu)
\Gamma(1+\nu+b)}
+ \mathcal{O}\left(\frac{1}{\nu^2}\right)
\label{jacobiphase:Mapprox2}
\end{dmath}
as $\nu \to \infty$.  
By inserting the asymptotic approximation
\begin{equation}
\frac{\Gamma(1+\nu+a)\Gamma(1+\nu+a+b)}{\Gamma(1+\nu)\Gamma(1+\nu+b)}
\sim p^{2a} \ \ \mbox{as} \  \ \nu\to\infty
\end{equation}
into  (\ref{jacobiphase:Mapprox2}), we obtain
\begin{equation}
\left(M^{(a,b)}(t,\nu)\right)^2
\sim
\frac{2}{\pi}
+ \mathcal{O}\left(\frac{1}{\nu^2}\right)
\ \ \mbox{as} \ \ \nu \to \infty.
\label{jacobiphase:Mapprox3}
\end{equation}
Once  again, (\ref{jacobiphase:Mapprox3}) 
 only holds for $t$ in the interior of the interval $(0,\pi)$
and $(a,b) \in (-1/2,1/2)$.
Figure~\ref{figure:phase} contains plots of the function
$\frac{\pi}{2} \left(M_\nu^{(a,b)}\right)^2$ 
for three sets of the parameters $a$, $b$ and $\nu$.
A straightforward, but somewhat tedious, computation
shows that the Wronskian of the pair $\tilde{P}_\nu^{(a,b)}$, $\tilde{Q}_\nu^{(a,b)}$ is
\begin{equation}
W_\nu^{(a,b)} = \frac{2p}{\pi}.
\end{equation}
From this, (\ref{jacobiphase:Mapprox3}) and (\ref{phase:mpsip}),
%
we conclude that for $t \in (\delta,\pi-\delta)$ and
$a,b \in (-1/2,1/2)$,
\begin{equation}
\psi^{(a,b)}(t,\nu) \sim pt + c + \mathcal{O}\left(\frac{1}{\nu^2}\right)
\ \ \mbox{as} \ \ \nu\to\infty
\end{equation}
with $c$ a constant.  Without loss of generality, we can
assume $c$ is in the interval $(0,2\pi)$.   This is Formula~(\ref{introduction:asym}).

\begin{figure}[t!!]
\begin{center}

\includegraphics[width=.32\textwidth]{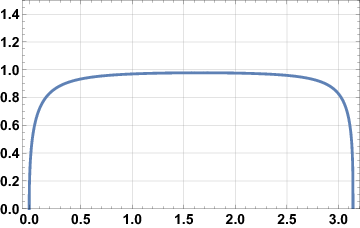}
\hfil
\includegraphics[width=.32\textwidth]{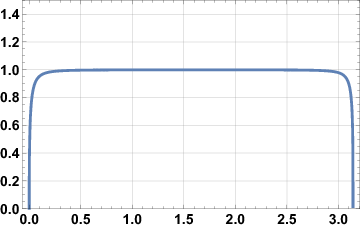}
\hfil
\includegraphics[width=.32\textwidth]{phase2.png}
\end{center}
\caption{On the left is a plot of the function $\frac{\pi}{2}\left(M_\nu^{(a,b)}\right)^2$
when $a=\frac{1}{4}$, $b=\frac{1}{3}$ and
$\nu = 1$.  In the middle is a plot of the same function
when $a=\frac{1}{4}$, $b=\frac{1}{3}$ and
$\nu = 10$.  On the right, is a plot of it
 when $a=\frac{1}{4}$, $b=\frac{1}{3}$ and $\nu = 100$.}
\label{figure:phase}
\end{figure}

\label{section:jacobiphase}
\end{subsection}

\begin{subsection}{Interpolative decompositions}

We say a factorization 
\begin{equation}
A = B R
\label{interp:1}
\end{equation}
of an $n\times m$ matrix $A$ is an interpolative decomposition 
 if $B$ consists of a collection of columns of the matrix $A$.
If (\ref{interp:1}) does not hold exactly, but instead
\begin{equation}
\|A - BR \|_2 \leq \epsilon \|A\|_2,
\end{equation}
then we say that (\ref{interp:1}) is an $\epsilon$-accurate
interpolative decomposition.    In our implementation of
the algorithms of this paper, a variant of the algorithm of 
\cite{Cheng-Gimbutas-Martinsson-Rokhlin} is used in order to form
interpolative decompositions.
\label{section:interp}
\end{subsection}

\label{section:preliminaries}
\end{section}

\begin{subsection}{Piecewise Chebyshev interpolation}

For our purposes, the $k$-point Chebyshev grid on the interval
$(\sigma_1,\sigma_2)$ is the set of points
\begin{equation}
\left\{
\frac{\sigma_2-\sigma_1}{2} \cos\left(\frac{\pi}{k-1} (k-i)\right) + \frac{\sigma_2+\sigma_1}{2}
: i=  1,\ldots,k \right\}.
\end{equation}
As is well known, given the values of a  polynomial $p$ of degree
$n-1$ at these nodes, its value at any point $t$ in the interval
$[\sigma_1,\sigma_2]$ can be evaluated in a numerically stable fashion
using the barcyentric Chebyshev interpolation formula
\begin{equation}
p(t) = \prod_{j=1}^{k}  \frac{w_j}{t-x_j} p(x_j)
\ \bigg / \ \ 
\prod_{j=1}^{k}  \frac{w_j}{t-x_j},
\label{cheby:bary}
\end{equation}
where $x_1,\ldots,x_k$ are the nodes of the $k$-point Chebyshev grid
on $(\sigma_1,\sigma_2)$ and 
\begin{equation}
w_j = \begin{cases}
1 & 1 < j < k \\
\frac{1}{2} & \mbox{otherwise}.
\end{cases}
\end{equation}
See, for instance, \cite{Trefethen}, for a detailed discussion
of Chebyshev interpolation and the numerical properties of Formula~(\ref{cheby:bary}).

We say that a function $p$ is 
a piecewise polynomial of degree $k-1$ on the collection of intervals
\begin{equation}
\left(\sigma_1,\sigma_2\right), \left(\sigma_2,\sigma_3\right), \ldots, \left(\sigma_{m-1},\sigma_m\right),
\label{cheby:intervals}
\end{equation}
where $\sigma_1 < \sigma_2 < \ldots < \sigma_m$,
provided its restriction to each $(\sigma_j,\sigma_{j+1})$ is a polynomial
of degree less than or equal to $k-1$.    We refer to the collection of 
points
\begin{equation}
\left\{
\frac{\sigma_{j+1}-\sigma_j}{2} \cos\left(\frac{\pi}{k-1} (k-i)\right) + 
\frac{\sigma_{j+1}+\sigma_j}{2}
: j=1,\ldots,m-1,\ i =1,\ldots,k \right\}
\label{cheby:disc}
\end{equation}
as the $k$-point piecewise Chebyshev grid on the intervals (\ref{cheby:intervals}).
Because the last point  of
each interval save $(\sigma_{m-1},\sigma_m)$ coincides with the first point on the succeeding interval,
there are $(k-1)(m-1) + 1$ such points.
Given the values of a piecewise polynomial $p$ of degree $k-1$ at these nodes,
its value at any point $t$ on the interval $[\sigma_1,\sigma_m]$ can be calculated
by identifying the interval containing $t$ and applying an appropriate modification
of Formula~(\ref{cheby:bary}).  Such a procedure requires at most 
$\mathcal{O}\left(m+k\right)$ operations, typically
requires $\mathcal{O}\left(\log(m)+k\right)$ when a binary search is used to identify
the interval containing $t$, and requires only 
 $\mathcal{O}(k)$ operations in the fairly common case that 
the  collection (\ref{cheby:intervals}) is of a form which allows the correct
interval to be identified immediately.  

We refer to the device of representing smooth functions
via their values at the  nodes (\ref{cheby:disc})
as the  piecewise Chebyshev discretization scheme of order $k$
on  (\ref{cheby:intervals}), or sometimes
just the piecewise Chebyshev discretization scheme  on (\ref{cheby:intervals})
when the value of $k$ is readily apparent.
Given such a scheme and  an ordered collection of points 
\begin{equation}
\sigma_1 \leq t_1 \leq t_2 \leq \ldots < t_N \leq \sigma_m,
\label{cheby:nodes}
\end{equation}
 we refer to the $N \times ((m-1)(k-1)+1)$ matrix
which maps the vector of values
\begin{equation}
\left(
\begin{array}{c}
p(x_{1,1}) \\
p(x_{2,1}) \\
\vdots    \\
p(x_{k,1}) \\
p(x_{2,2}) \\
\vdots    \\
p(x_{k,2}) \\
\vdots    \\
p(x_{k,m-1}) \\
p(x_{2,m}) \\
\vdots    \\
p(x_{n,m}) \\
\end{array}
\right),
 \end{equation}
where $p$ is any piecewise polynomial of degree $k-1$
on the interval (\ref{cheby:intervals}) and 
\begin{equation}
x_{i,j} = 
\frac{\sigma_{j+1}-\sigma_j}{2} \cos\left(\frac{\pi}{k-1} (k-i)\right) + 
\frac{\sigma_{j+1}+\sigma_j}{2},
\end{equation}
to the vector of values
\begin{equation}
\left(
\begin{array}{c}
p(t_1) \\
p(t_2) \\
\vdots    \\
p(t_N)
\end{array}
\right)
\end{equation}
as the matrix interpolating 
piecewise polynomials of degree $k-1$ on
(\ref{cheby:intervals}) to the nodes (\ref{cheby:nodes}).
This matrix is block diagonal; in particular, it is of the form
\begin{equation}
\left(
\begin{array}{cccc}
B_1 & 0 & 0 & 0 \\
0  & B_2 & 0 & 0 \\
0 &  0 & \ddots & 0 \\
0 &  0 & 0 & B_{m-1} \\
\end{array}
\right)
\end{equation}
with $B_j$ corresponding to the interval
 $(\sigma_j,\sigma_{j+1})$.  The dimensions of $B_j$ are
$n_j \times k$, where $n_j$ is the number of points from
(\ref{cheby:nodes}) which lie in $\left[\sigma_j,\sigma_{j+1}\right]$.
The entries of each block can be calculated from (\ref{cheby:bary}),
with the consequence that this matrix can be applied in 
 $\mathcal{O}\left(nk\right)$ operations without explicitly forming it.

\label{section:cheby}
\end{subsection}

\begin{subsection}{Bivariate Chebyshev interpolation}

Suppose that
\begin{equation}
x_1,\ldots,x_{M_1},
\label{bivariate:xs}
\end{equation}
where $M_1=(m-1)(k_1-1)+1$, is the $k_1$-point piecewise Chebyshev grid on the collection of intervals
\begin{equation}
\left(\sigma_1,\sigma_2\right), \left(\sigma_2,\sigma_3\right), \ldots, \left(\sigma_{m_1-1},\sigma_{m_1}\right),
\label{bivariate:ints1}
\end{equation}
and that 
\begin{equation}
y_1,\ldots,y_{M_2},
\label{bivariate:ys}
\end{equation}
where $M_2=(m_2-1)(k_2-1)+1$,
is the $k_2$-point piecewise Chebyshev grid on the collection of intervals
\begin{equation}
\left(\zeta_1,\zeta_2\right), \left(\zeta_2,\zeta_3\right), \ldots, \left(\zeta_{m_2-1},\zeta_{m_2}\right).
\label{bivariate:ints2}
\end{equation}
Then we refer to the device of representing
smooth bivariate functions $f(x,y)$ defined on the rectangle
$(\sigma_1,\sigma_{m_1}) \times (\zeta_1,\zeta_{m_2})$ 
via their values at the tensor product 
\begin{equation}
\left\{
\left(x_i,y_j\right) : i=1,\ldots,M_1, \ j=1,\ldots,M_2
\right\}
\label{bivariate:tensor}
\end{equation}
of the piecewise  Chebyshev grids (\ref{bivariate:xs})
and (\ref{bivariate:ys}) as a bivariate Chebyshev discretization
scheme.  Given the values of $f$ at the discretization nodes (\ref{bivariate:tensor}),
its value at any point $(x,y)$ in the rectangle 
$(\sigma_1,\sigma_{m_1}) \times (\zeta_1,\zeta_{m_2})$ 
can be approximated through repeated applications of the barcyentric Chebyshev interpolation
formula.  The cost of such a procedure is $\mathcal{O}\left(k_1k_2\right)$, once
the rectangle $\left(\sigma_i,\sigma_{i+1}\right) \times
\left(\zeta_j,\zeta_{j+1}\right)$ containing $(x,y)$ has been located.
The cost of locating the correct rectangle is $\mathcal{O}\left(m_1+m_2\right)$
in the worst case, $\mathcal{O}\left(\log(m)_1+\log(m_2)\right)$ in typical cases,
and $\mathcal{O}\left(1\right)$ when the particular forms of the collections
(\ref{bivariate:ints1}) and (\ref{bivariate:ints2}) allow for the
immediate  identification of the correct rectangle.

\label{section:bivariate}
\end{subsection}

\begin{section}{Computation of the Phase and Amplitude Functions} 

Here, we describe a method for 
calculating representations of the nonoscillatory phase and amplitude functions 
$\psi^{(a,b)}(t,\nu)$
and $M^{(a,b)}(t,\nu)$.
This procedure takes as input the parameters $a$ and $b$ as well as an integer
 $N_{\mbox{\tiny max}} > 27$. 
The representations of $M^{(a,b)}$ and $\psi^{(a,b)}$
allow for their evaluation for all values of
$\nu$ such that
\begin{equation}
27 \leq \nu \leq N_{\mbox{\tiny max}}
\label{precomp:nurange}
\end{equation}
and all arguments $t$ such that
\begin{equation}
\frac{1}{N_{\mbox{\tiny max}}} \leq t \leq \pi - \frac{1}{N_{\mbox{\tiny max}}}.
\label{precomp:trange}
\end{equation}
The smallest and largest zeros of $\tilde{P}_\nu^{(a,b)}$ on the interval
$(0,\pi)$ are greater
than $\frac{1}{\nu}$ and less than $\pi -  \frac{1}{\nu}$, respectively (see, for instance,
Chapter~18 of \cite{DLMF} and its references).    In particular, this means that the range
(\ref{precomp:nurange}) suffices to evaluate $\tilde{P}_\nu^{(a,b)}$ at the values
of $t$ corresponding to the nodes of the
$N_{\mbox{\tiny max}}$-point Gauss-Jacobi quadrature rule.  
If necessary, this range could be expanded or various asymptotic expansion
could be used to evaluate $\tilde{P}^{(a,b)}(t,\nu)$ for values of $t$ outside of it.

The phase and amplitude functions are represented via their values
at the nodes of a tensor product of piecewise Chebyshev grids.
To define them, we first let $m_\nu$ be the least integer such that
\begin{equation}
3^{m_\nu} \geq N_{\mbox{\tiny max}}
\end{equation}
and let $m_t$ be equal to twice the least integer $l$ such that
\begin{equation}
\frac{\pi}{2} 2^{-l+1} \leq \frac{1}{N_{\mbox{\tiny max}}}.
\end{equation}
Next, we define  $\beta_j$ for $j=1,\ldots,m_\nu$ by
\begin{equation}
\beta_j = \max\left\{ 3^{j+2}, N_{\mbox{\tiny max}}\right\},
\end{equation}
$\alpha_i$ for $i=1,\ldots,m_t/2$ via
\begin{equation}
\alpha_i = \frac{\pi}{2}  2^{i-m_t/2}
\end{equation}
and $\alpha_i$ for $i=m_t/2+1, \ldots,m_t$ via
\begin{equation}
\beta_i = \pi - \frac{\pi}{2}  2^{m_t/2+1-i}.
\end{equation}
We note that the points $\{\beta_i\}$ cluster near  $0$ and $\pi$,
where the phase and amplitude functions are singular.
Now we let 
\begin{equation}
\tau_1,\ldots,\tau_{M_t}
\label{precomp:tnodes}
\end{equation}
denote $16$-point piecewise Chebyshev gird on the intervals
\begin{equation}
\left(\alpha_1,\alpha_2\right), \left(\alpha_2,\alpha_3\right), \ldots,
\left(\alpha_{m_t-1},\alpha_{m_t}\right).
\label{precomp:tints}
\end{equation}
There are
\begin{equation}
M_t = 15( m_t-1) + 1
\end{equation}
such points.  In a similar fashion, we let 
\begin{equation}
\gamma_1,\ldots,\gamma_{M_\nu}
\label{precomp:nunodes}
\end{equation}
denote the nodes of the $24$-point piecewise Chebyshev grid
on the intervals
\begin{equation}
\left(\beta_1,\beta_2\right), \left(\beta_2,\beta_3\right), \ldots,
\left(\beta_{{m_\nu}-1},\beta_{m_{\nu}}\right).
\label{precomp:nuints}
\end{equation}
There are 
\begin{equation}
M_t = 23( m_\nu-1) + 1
\end{equation}
such points.

For each node $\gamma$ in the set
(\ref{precomp:nunodes}),  we solve the ordinary differential equation (\ref{phase:Neq})
with $q$ taken to be the coefficient (\ref{jacobieq:coef}) for Jacobi's modified differential
equation using   a variant
of the integral equation method of \cite{Greengard}.  We note, though, that
any standard method capable of coping with stiff problems should suffice.
We determine a unique solution   by specifying
 the values of 
\begin{equation}
\left(M^{(a,b)}(t,\gamma)\right)^2
\label{precomp:msq}
\end{equation}
 and its first two derivatives
at a point on the interval (\ref{precomp:trange}).  These values
are calculated using an asymptotic expansion for (\ref{precomp:msq})
derived from 
(\ref{asym:p}) and (\ref{asym:q}) when $\gamma$ is large,
and using series expansions for $P_\gamma^{(a,b)}$ and $Q_\gamma^{(a,b)}$
when $\gamma$ is small.  We refer the interested reader
to our code (which we have made publicly available) for details.
Solving (\ref{phase:Neq})  gives us the values of the function
for each $t$ in the set (\ref{precomp:tnodes}).  To obtain the
values of $\psi^{(a,b)}(t,\gamma)$ for each $t$ in the set (\ref{precomp:tnodes}) , we first calculate
the values of
\begin{equation}
\frac{d}{dt} \psi^{(a,b)}(t,\gamma)
\end{equation}
via (\ref{phase:mpsip}).      Next, we obtain
the values of the function $\tilde{\psi}$  defined via
\begin{equation}
\tilde{\psi}(t) = \int_{\alpha_1}^t 
\frac{d}{ds} \psi^{(a,b)}(s,\gamma)\ ds
\end{equation}
at the nodes (\ref{precomp:tnodes}) via spectral integration.  There is
an as yet unknown constant connecting $\tilde{\psi}$ with 
 $\psi^{(a,b)}(t,\gamma)$; that is,
\begin{equation}
\psi^{(a,b)}(t,\gamma) = \tilde{\psi}(t) + C.
\end{equation}
To evaluate $C$, we first use a combination of asymptotic and series
expansions to calculate $\tilde{P}_\nu^{(a,b)}$ at the point $\alpha_1$.
Since $\widetilde{\psi}(\alpha_1) =0$, it follows that 
\begin{equation}
\tilde{P}_\gamma^{(a,b)}(\alpha_1) = M^{(a,b)}(\alpha_1,\gamma) \cos( C),
\label{precomp:c}
\end{equation}
and $C$ can be calculated in the obvious fashion.  Again, we refer
the interested reader to our source code for the details
of the asymptotic expansions we use to evaluate 
$\tilde{P}_\nu^{(a,b)}(\alpha_1)$.

The ordinary differential equation (\ref{phase:Neq}) is solved for 
$\mathcal{O}\left(\log\left(N_{\mbox{\tiny max}}\right)\right)$ 
values of $\nu$, and 
the phase and amplitude functions are calculated at 
$\mathcal{O}\left(\log\left(N_{\mbox{\tiny max}}\right)\right)$ 
values of $t$ each time it is solved.    This makes the total running time
of the procedure just described
\begin{equation}
\mathcal{O}\left(\log^2\left(N_{\mbox{\tiny max}}\right)\right).
\end{equation}
Once the values of $\psi_{\nu}^{(a,b)}$ and $M_\nu^{(a,b)}$
have been calculated at the tensor product
of the piecewise Chebyshev grids (\ref{precomp:tnodes}) and (\ref{precomp:nunodes}),
they can be evaluated for any $t$ and $\nu$ via repeated application
of the barycentric
Chebyshev interpolation formula  in
a number of operations which
is independent of $\nu$ and $t$.  The rectangle in the discretization scheme
which contains the point $(t,\nu)$ can be determined in $\mathcal{O}\left(1\right)$
operations
owing to the special form of (\ref{precomp:tints}) and (\ref{precomp:nuints}).
The value of $\tilde{P}_\nu^{(a,b)}$
can then be calculated via Formula~(\ref{jacobiphase:p}),
also in a number of operations which is independent of $\nu$ and $t$.

\begin{remark}
Our  choice of the particular bivariate piecewise Chebyshev discretization
scheme  used to represent the functions
$\psi^{(a,b)}$ and  $M^{(a,b)}$ was informed by extensive numerical
experiments.  However, this does not preclude
the possibility that there are other similar schemes which provide better levels of accuracy
and efficiency.     Improved mechanisms for the representation
of phase functions are currently being  investigated by the authors.
\end{remark}


\label{section:precomp}
\end{section}

\begin{section}{Computation of Gauss-Jacobi Quadrature Rules}

In this section, we describe our algorithm for the calculation
of Gauss-Jacobi quadrature rules (\ref{introduction:quadrule}),
and for the modified Gauss-Jacobi quadrature rules
 (\ref{introduction:modrule}).
It takes as input the length $n$
of the desired quadrature
rule as well as  parameters $a$ and $b$ in the interval $\left(-\frac{1}{2},\frac{1}{2}\right)$
and returns the nodes $x_1,\ldots,x_n$ and weights $w_1,\ldots,w_n$
of (\ref{introduction:quadrule}).    

For these computations, we do not rely on the 
expansions of $\psi^{(a,b)}$ and $M^{(a,b)}$ produced using
the procedure of Section~\ref{section:precomp}.
Instead, as a first step,
 we calculate the values of 
$\psi^{(a,b)}(t,n)$  and $M^{(a,b)}(t,n)$ for each $t$  in the set (\ref{precomp:tnodes}) by solving
the ordinary differential equation (\ref{phase:Neq}).  The procedure
  is the same as that used in the algorithm of the preceding section.  
Because the degree $\nu = n$ of the phase function used by the algorithm
of this section is fixed, we break from our notational convention 
and use $\psi^{(a,b)}_n$ to denote the relevant phase function; that is,
\begin{equation}
\psi^{(a,b)}_n(t) = \psi^{(a,b)}\left(t,n\right).
\end{equation}
  We next calculate the inverse function 
$\left(\psi_n^{(a,b)}\right)^{-1}$ 
of the phase function $\psi_n^{(a,b)}$ as follows.
First, we define $\omega_i$ for $i=1,\ldots,m_t$ via the formula
\begin{equation}
\omega_i = \psi_n^{(a,b)}\left(\alpha_i\right),
\end{equation}
where $m_t$ and $\alpha_1,\ldots,\alpha_{m_t}$ are as in the preceding section.
Next, we form the collection
\begin{equation}
\xi_1,\ldots,\xi_{M_t}
\label{quad:nodes}
\end{equation}
of points obtained by taking the
union of $16$-point Chebyshev grids on each of the intervals
\begin{equation}
\left(\omega_1,\omega_2\right),
\left(\omega_2,\omega_3\right), \ldots,
\left(\omega_{m_t-1},\omega_{m_t}\right).
\end{equation}
Using Newton's method, we compute the value of the inverse function
of $\psi_n^{(a,b)}$ at each of the points (\ref{quad:nodes}).  
More explicitly, we traverse the set (\ref{quad:nodes}) in descending order,
and for each node $\xi_j$ solve the equation
\begin{equation}
\psi^{(a,b)}_n(y_j) = \xi_j
\end{equation}
for $y_j$ so as to determine the value of 
\begin{equation}
\left(\psi_n^{(a,b)}\right)^{-1}(\xi_j).
\end{equation}
The values of the derivative of $\psi_n^{(a,b)}$ are a by-product
of the procedure used to construct $\psi_n^{(a,b)}$, and so the evaluation
of the derivative of the phase function is straightforward.
Once the values of the inverse function at the nodes (\ref{quad:nodes})
are known,  barycentric Chebyshev interpolation enables us to 
 evaluate $\psi^{-1}$ at any point
on the interval $\left(\omega_1,\omega_{m_t}\right)$.
From (\ref{jacobiphase:p}), we see that the $k^{th}$ zero $t_k$ of $\tilde{P}_n^{(a,b)}$
is given by  the formula
\begin{equation}
t_k = \left(\psi_n^{(a,b)}\right)^{-1}\left(\frac{\pi}{2} + k \pi\right).
\label{quad:roots}
\end{equation}
These are the nodes of the modified Gauss-Jacobi quadrature rule,
and we  use (\ref{quad:roots}) to compute them.
 The corresponding weights are given by the formula
\begin{equation}
w_k = \frac{\pi}{ \frac{d}{dt} \psi_n^{(a,b)}(t_k) }.
\label{quad:wht1}
\end{equation}
The nodes of the Gauss-Jacobi quadrature rule are given by
\begin{equation}
x_j = \cos\left(t_{n-k+1}\right), 
\end{equation}
and the corresponding weights  are given by
\begin{equation}
w_k = \frac{\pi 2^{a+b+1}  \sin\left(\frac{t}{2}\right)^{2a+1} \cos\left(\frac{t}{2}\right)^{2b+1}}
{\frac{d}{dt} \psi_n^{(a,b)}(t_{n-k+1}) }.
\label{quad:wht2}
\end{equation}
Formula (\ref{quad:wht2}) can be obtained by combining
 (\ref{jacobiphase:p}) with Formula~(15.3.1)
in Chapter~15 of \cite{Szego}.
That
(\ref{quad:wht1}) gives the values of the weights for (\ref{introduction:quadrule})
is then obvious from a comparison of the form of that rule
and (\ref{introduction:modrule}).

The cost of computing the phase function and its inverse is $\mathcal{O}\left(\log(n)\right)$,
while the cost of computing each node and weight pair is independent of $n$.  
Thus the algorithm of this section
has an asymptotic running time of $\mathcal{O}\left(n\right)$.

\begin{remark}
The phase function $\psi_n^{(a,b)}$ is a monotonically increasing function
whose condition number of evaluation is small, whereas 
 the Jacobi polynomial $P_n^{(a,b)}$ it represents is highly oscillatory
with  a large condition number of  evaluation.  
As a consequence, many numerical computations, including 
the numerical calculation of the zeros of $P_n^{(a,b)}$, benefit from a change to phase
space.  See \cite{BremerZeros} for a more thorough discussion of this observation.
\end{remark}

\label{section:quad}
\end{section}

\begin{section}{Application of the Jacobi Transform and its Inverse}

In this section, we describe a procedure for applying 
the $n^{th}$ order Jacobi transform --- which is represented by
the $n\times n$ matrix
 $\mathcal{J}_n^{(a,b)}$ defined via
(\ref{introduction:jn}) ---  to a vector $x$.
 Since 
 $\mathcal{J}_n^{(a,b)}$ is orthogonal,  the inverse Jacobi transform simply consists of 
applying its transpose.  We omit a discussion of the minor and obvious
changes to the algorithm of this section required to do so.

A precomputation, which includes as a step the  procedure for the construction
of the  nonoscillatory phase and amplitude functions
described Section~\ref{section:precomp},
 is necessary before our  algorithm for applying the forward Jacobi transform
can be used.  We first discuss
the outputs of the precomputation procedure 
and our method for using them to apply the Jacobi transform.
Afterward, we describe the precomputation procedure in detail.

The precomputation phase of this algorithm
takes as input the order $n$ of the transform
and parameters $a$ and $b$ in the interval $\left(-\frac{1}{2},\frac{1}{2}\right)$.
Among other things, it produces the nodes $t_1,\ldots,t_n$ and 
weights $w_1,\ldots,w_n$ of the  $n$-point modified Gauss-Jacobi quadrature rule,
and a second collection of points
$x_1,\ldots,x_n$ such that for each $j=1,\ldots,n$, $x_j$ is the node
from the equispaced grid
\begin{equation}
0,\ 2\pi \cdot \frac{1}{n},\ 2\pi\cdot \frac{2}{n},\ \ldots,\ 2\pi\cdot\frac{n-1}{n}
\label{transform:equispaced}
\end{equation}
which is  closest to   $t_j$.  In addition, the precomputation phase constructs
a complex-valued $n \times r$ matrix $L$, a
complex-valued $r \times n$ matrix $R$ and  a real-valued $n \times 26$
matrix $V$ such that
\begin{equation}\label{eqn:FLR}
\mathcal{J}_n^{(a,b)} \approx 
\left(
\real\left(
\mathcal{F}_n\otimes (L  R)
\right)
+
\left(\begin{array}{cc} V & 0\end{array}\right)
\right)
W,
\end{equation}
where  $r$ is the numerical rank of a matrix we define shortly,
$\mathcal{F}_n$ is the $n \times n$ matrix whose $(j,k)$ entry is
\begin{equation}
 \exp( i x_j (k-1)),
\end{equation}
$W$ is the $n\times n$ matrix
\begin{equation}
W = \left(
\begin{array}{cccc}
\sqrt{w_1} & 0          & 0      & 0           \\
0          & \sqrt{w_2} & 0      & 0           \\
0          & 0          & \ddots & 0           \\
0          & 0          & 0      & \sqrt{w_n} \\
\end{array}
\right),
\end{equation}
and $\otimes$ denotes the Hadamard product.
  The matrix $\mathcal{F}_n$ 
can be applied to a vector in
 $\mathcal{O}\left(n \log(n)\right)$ operations by executing $r$ fast Fourier transforms.
Moreover, if we use $u_1,\ldots,u_r$ to denote the $n \times 1$ 
matrices which comprise  the columns of
the matrix $L$ and $v_1,\ldots,v_r$ to denote the $1 \times n$ matrices
which comprise the rows of the matrix $R$, then
\begin{equation}
LR = u_1v_1 + \ldots + u_r v_r
\end{equation}
and 
\begin{equation}
\mathcal{F}_n \otimes \left(LR\right)
= 
D_{v_1} \mathcal{F}_n D_{u_1} + \cdots + 
D_{v_r} \mathcal{F}_n D_{u_r},
\label{transform:decomp2}
\end{equation}
where  $D_u$ denotes the $n \times n$ diagonal matrix with entries $u$ on the diagonal.
From (\ref{transform:decomp2}), we see that 
$\mathcal{F}_n \otimes \left(LR\right)$ can be applied in $\mathcal{O}\left(rn\log(n)\right)$
operations using the Fast Fourier transform.  This 
is the approach used in \cite{TownsendNUFFT} to 
rapidly apply discrete nonuniform Fourier transforms.
Since the matrices $\left(\begin{array}{cc} V & 0\end{array}\right)$
and $W$
can each be applied to a vector in $\mathcal{O}\left(n\right)$ operations,
the asymptotic running time of this procedure for applying 
the forward Jacobi transform  is
$\mathcal{O}\left(r n\log(n) \right)$.

We now describe the precomputation step.
First, the procedure  of Section~\ref{section:precomp}
is used to construct
 nonoscillatory phase  and amplitude functions $\psi^{(a,b)}$ and $M^{(a,b)}$;
the parameter   $N_{\mbox{\tiny max}}$ is taken to be  $n$.
In what follows, we let 
 $M_t$, $M_\nu$, $\alpha_j$, $\beta_j$, $\tau_j$ and $\gamma_k$ be as in Section~\ref{section:precomp}.
Next, the procedure of Section~\ref{section:quad} is used to construct
the nodes and weights of the modified Gauss-Jacobi quadrature rule
of order $n$.     We now form an interpolative decomposition
\begin{equation}
\left(
\begin{array}{c}
A^{(1)}\\
A^{(2)}
\end{array}
\right)
\approx
\left(
\begin{array}{c}
B^{(1)}\\
B^{(2)}
\end{array}
\right) \widetilde{R}
\label{transform:decomp}
\end{equation}
(see Section~\ref{section:interp}) of the matrix

\begin{equation}
\left(
\begin{array}{c}
A^{(1)}\\
A^{(2)}
\end{array}
\right),
\label{transform:matrix}
\end{equation}
where $A^{(1)}$ is the $M_t \times M_\nu$ matrix  whose  $(j,k)$ entry is given by
\begin{equation}
A^{(1)}_{jk} = 
M^{(a,b)}(\tau_j,k-1)
\exp\left(i \left( \psi^{(a,b)}(\tau_j,\gamma_k) -  \tau_j\gamma_k \right)\right),
\end{equation}
and $A^{(2)}$ is  the $16 \times M_\nu$ matrix  whose entries are 
\begin{equation}
A^{(2)}_{jk} = 
\exp\left(  i\ \sigma_j\  \frac{\gamma_k}{N_{\mbox{\tiny max}}}
\right)
\end{equation}
with $\sigma_1\ldots,\sigma_{16}$ are the nodes of the  $16$-point Chebyshev grid
on the interval $(0,2\pi)$.
We use the
procedure of \cite{Cheng-Gimbutas-Martinsson-Rokhlin}
with  $\epsilon = 10^{-12}$ requested accuracy to do so.
We view this procedure as choosing a collection
\begin{equation}
\gamma_{s_1} < \ldots <\gamma_{s_r}
\label{precomp:ranknodes}
\end{equation}
of the nodes from the set 
(\ref{precomp:nunodes}), where $r$ is the numerical rank of  matrix
(\ref{transform:matrix}) to precision $\epsilon$,
 as determined by the interpolative  decomposition procedure,
and constructing a device (the matrix $\widetilde{R}$) for evaluating functions of the form
\begin{equation}
f(\nu) =  M^{(a,b)}(t,\nu) \exp\left( i \left(\psi^{(a,b)}(t,\nu) - t\nu\right) \right)
\label{precomp:fun1}
\end{equation}
and
\begin{equation}
g(\nu) =  \exp\left( i t \frac{\nu}{N_{\mbox{\tiny max}}} \right),
\label{precomp:fun2}
\end{equation}
where $t$ is any point in the interval (\ref{precomp:trange}),
at the nodes (\ref{precomp:nunodes}) 
 given their values at the nodes  (\ref{precomp:ranknodes}).
That $t$ can take on any value in the 
interval (\ref{precomp:trange}) and not just the special
values (\ref{precomp:tnodes}) is a consequence of the fact
that the piecewise bivariate Chebyshev discretization scheme described
in the preceding section suffices to represent these functions.
We expect (\ref{transform:matrix}) to be low rank from
the discussion in the introduction and the estimates
of Section~\ref{section:jacobiphase}.  We factor
the augmented matrix (\ref{transform:matrix}) rather than
$A^{(1)}$ because in later steps of the precomputation
procedure, we use
the matrix $\widetilde{R}$ to 
interpolate functions of the form (\ref{precomp:fun2})
as well as those of the form (\ref{precomp:fun1}).

In what follows, we denote by $I_{\mbox{\tiny left}}$  the $n \times M_t$ matrix which
interpolates piecewise polynomials of degree $15$ on the intervals
(\ref{precomp:tints}) to their values
at the points $t_1,\ldots,t_n$.  
See Section~\ref{section:cheby} for a careful definition of the 
 matrix $I_{\mbox{\tiny left}}$ and a discussion of its properties.
Similarly, we let 
$I_{\mbox{\tiny right}}$ be the $M_\nu \times n$ matrix
\begin{equation}
I_{\mbox{\tiny right}} = \left(
\begin{array}{cc}
0 & I_{\mbox{\tiny right}}^{(1)}
\end{array}
\right),
\end{equation}
where $I_{\mbox{\tiny right}}^{(1)}$ is the 
 $M_{\nu} \times (n-27)$ matrix which interpolates piecewise
polynomials of degree $23$ on the intervals (\ref{precomp:nuints})
 to their values  at the points
$27,28,29,\ldots,n-1$
(recall that the phase and amplitude functions represent the Jacobi polynomials
of degrees greater than or equal to $27$).    
We next form the $M_t \times r$ matrix whose $(j,k)$ entry is 
\begin{equation}
M^{(a,b)}(\tau_j,\gamma_{s_k})
\exp\left(
i\left(
\psi\left(\tau_j,\gamma_{s_k}\right) - \tau_j \gamma_{s_k}\right)
\right)
\end{equation}
and multiply it on the left by $I_{\mbox{\tiny left}}$ so as to form
the $n \times r$ matrix whose $(j,k)$ entry is
\begin{equation}
M^{(a,b)}(t_j,\gamma_{s_k})
\exp\left(
i\left(
\psi\left(t_j,\gamma_{s_k}\right) - t_j \gamma_{s_k}\right)
\right).
\end{equation}
For $j=1,\ldots,n$ and $k=1,\ldots,r$, we scale the 
$(j,k)$ entry of this matrix by
\begin{equation}
\exp\left(
i 
(t_j-x_j) \gamma_{s_k}
\right)
\label{transform:fun3}
\end{equation}
in order to form the  $n \times r$ matrix $L$ whose $(j,k)$ entry is 
\begin{equation}
L_{jk} = M^{(a,b)}(t_j,\gamma_{s_k})
\exp\left(
i\left(
\psi\left(t_j,\gamma_{s_k}\right) - x_j \gamma_{s_k}\right)
\right).
\end{equation}
We multiply the $r \times M_\nu$ matrix $\widetilde{R}$ on the right
by $I_{\mbox{\tiny right}}$  to form the $r \times n$ matrix 
$R$.     
It is because of the scaling by (\ref{transform:fun3})
that we factor the augmented
matrix (\ref{transform:matrix}) rather than $A^{(1)}$.
The values of (\ref{transform:fun3}) can be  correctly
interpolated by the matrix $\widetilde{R}$ because
of the addition of $A^{(2)}$.

As the final step in our precomputation procedure,
we use the well-known three-term recurrence relations
satisfied by the Jacobi polynomials to form the 
$n \times k$ matrix
\begin{equation}
V = \left(
\begin{array}{ccccccc}
\tilde{P}_0^{(a,b)} (t_1) & \tilde{P}_1^{(a,b)} (t_1) & \cdots & \tilde{P}_{26}^{(a,b)} (t_1) \\
\tilde{P}_0^{(a,b)} (t_2) & \tilde{P}_1^{(a,b)} (t_2) & \cdots & \tilde{P}_{26}^{(a,b)} (t_2) \\
\vdots & \vdots & \ddots & \vdots \\
\tilde{P}_0^{(a,b)} (t_n) & \tilde{P}_1^{(a,b)} (t_n) & \cdots & \tilde{P}_{26}^{(a,b)} (t_n) \\
\end{array}
\right).
\end{equation}

The highest order term in the asymptotic running time of the precomputation
algorithm is $\mathcal{O}(rn)$, where $r$ is the rank of (\ref{transform:matrix}),
and it comes from the application of the interpolation
matrices $I_{\mbox{\tiny right}}$ and $I_{\mbox{\tiny left}}$.  Based on our
numerical experiments, we conjecture that
\begin{equation}
r = \mathcal{O}\left( \frac{\log(n)}{\log\log(n)}\right),
\end{equation}
in which case the asymptotic running time of the procedure is
\begin{equation}
\mathcal{O}\left(\frac{n\log^2(n)}{\log\log(n)}\right).
\end{equation}

\begin{remark}
There are 
many mechanisms for accelerating the procedure used here for the calculation
of interpolative decomposition (particularly, randomized
algorithms \cite{RandomSample} 
and methods based on adaptive cross approximation).
However, in our numerical experiments we  found that the cost of
the construction of the phase and amplitude functions dominated the running
time of our precomputation step for small $n$, and that
the cost of interpolation  did so for large $n$. 
As a consequence,  the optimization of 
our approach to forming interpolative decompositions is a low priority.
\end{remark}

\begin{remark}
 The rank of the low-rank approximation by interpolative decompositions is not optimally small. A fast truncated SVD based on the interpolative decomposition by Chebyshev interpolation can provide the optimal rank of the low-rank approximation in \eqref{eqn:FLR} and 
thereby speed up
 the calculation in \eqref{transform:decomp2} by a small factor
(see \cite{Li-Yang} for a comparison). However, we do not explore this in this paper.
\end{remark}

\begin{remark}
The procedure of this section can modified rather easily to apply
several transforms related to the Jacobi transform.  For instance,
the mapping which take the coefficients (\ref{introduction:transin})
in the expansion (\ref{introduction:transexp})
to values of $f$ at the nodes of a Chebyshev 
quadrature  can be applied using the procedure of this
section by simply modifying $t_1,\ldots,t_n$.  The Jacobi-Chebyshev
transform can then be computed from the resulting value using
the fast Fourier transform.
\end{remark}

\label{section:transform}
\end{section}


\begin{section}{Numerical Experiments}

In this section, we describe numerical experiments conducted to evaluate
the performance of the algorithms of this paper.
Our code was written in Fortran and  compiled with the GNU
Fortran compiler version 4.8.4.  It uses version 3.3.2 of the FFTW3 library \cite{FFTW05}
to apply the discrete Fourier transform. 
We have made our code, including 
that for all  of the experiments described here, available on GitHub at the following
address:
\begin{center}
\url{https://gitlab.com/FastOrthPolynomial/Jacobi.git}
\end{center}
%
All calculations were carried out 
on a single core of workstation computer equipped with an Intel Xeon E5-2697 processor 
running at 2.6 GHz and 512GB of RAM.



\begin{subsection}{Numerical evaluation of Jacobi polynomials}

In a first set of experiments, we measured the speed of the procedure
of Section~\ref{section:precomp} for the construction of nonoscillatory
phase and amplitude functions, and the speed and accuracy with which
Jacobi polynomials can be evaluated using them.

We did so in part by comparing  our approach to a reference method which combines 
the Liouville-Green type  expansions of Barratella and Gatteschi (see Section~\ref{section:gatteschi})
with the trigonometric expansions of Hahn (see Section~\ref{section:hahn}).
Modulo minor implementation details, the reference method
we used is the same as that used in \cite{Hale-Townsend} to  evaluate Jacobi polynomials and
is quite  similar to the approach of \cite{Bogaert-Michiels-Fostier} for the 
evaluation of Legendre polynomials.
To be more specific, we used the expansion (\ref{asym:p}) to evaluate
$\tilde{P}_\nu^{(a,b)}$ for all $t \geq 0.2$ and the trigonometric
expansion (\ref{hahn:p}) to evaluate it when $0 < t < 0.2$.  
The expansion (\ref{hahn:p}) is much more expensive than (\ref{asym:p});
 however, the use of
(\ref{asym:p}) is complicated by the fact that the higher order coefficients
are not known explicitly.  We calculated $16^{th}$ order Taylor expansion for the coefficients
$B_1$, $A_1$, $B_2$ and $A_2$, but the use of these approximations
limited the range of $t$ for we which we could apply (\ref{asym:p}).

\begin{figure}[b!!]
\begin{center}
\hfil
\includegraphics[width=.49\textwidth]{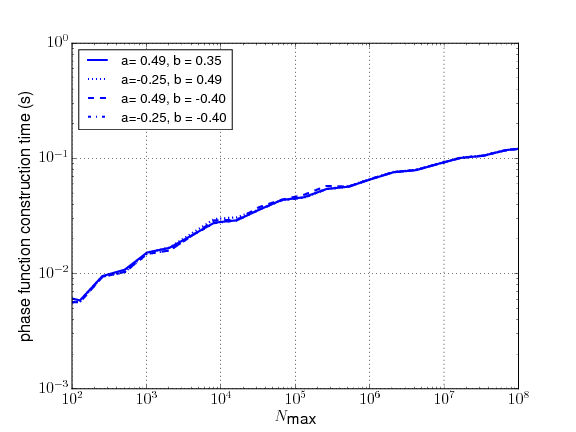}
\hfil 
\includegraphics[width=.49\textwidth]{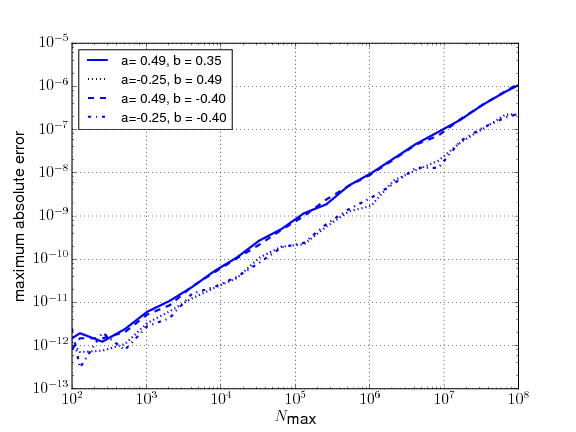}
\hfil
\end{center}
\label{figure:asym}
\caption{The plot on the left shows how the
 time required to construct the phase and amplitude functions varies with
$N_{\mbox{\tiny max}}$
for various values of the parameters $a$ and $b$.
The plot on the right shows the
maximum absolute errors which were observed when evaluating
Jacobi polynomials
as a function of  $N_{\mbox{\tiny max}}$ for various values of the parameters $a$ and $b$.}
\end{figure}

In order to perform these comparisons, 
we fixed $a = -1/4$ and $b=1/3$ and choose several values
of $N_{\mbox{\tiny max}}$.  For each such value,
we carried out the procedure 
for the construction of the phase function described in Section~\ref{section:precomp}
and evaluated $\tilde{P}_\nu^{(a,b)}$
for $200$ randomly chosen values of $t$ 
in the interval $(0,\pi)$
and $200$ randomly chosen values of $\nu$ in the interval $(1,N_{\mbox{\tiny max}})$
using both our algorithm and the reference method.
Table~\ref{table:asym} reports the results.  
For several different pairs of the parameters $a,b$, we repeated these experiments,
measuring the time required to construct the phase function and the 
maximum observed absolute errors.  Figure~\ref{figure:asym} displays the results.
We note that the condition number of evaluation of the Jacobi polynomials
increases with order, with the consequence that the obtained accuracy
of both the approach of this paper and the 
reference method decrease as a function of degree.
The errors observed in Table~\ref{table:asym} and Figure~\ref{figure:asym}
are consistent with this fact.  

\begin{table}[t!]
\small
\begin{center}
 \begin{tabular}{lccccc}
 \toprule
 $N_{\mbox{\tiny max}}$ & Phase function  & Avg. eval time & Avg. eval time& Largest absolute & Expansion size\\
 &construction time & algorithm & asymptotic &   error & (MB)\\
 &  & of this paper & expansions  &   &\\
 \midrule
100 & 1.15\e{-02} & 2.82\e{-07} & 3.46\e{-05} & 1.31\e{-12} & 1.90\e{-01}  \\
128 & 5.77\e{-03} & 7.82\e{-07} & 2.56\e{-05} & 8.89\e{-13} & 1.90\e{-01}  \\
256 & 9.46\e{-03} & 5.81\e{-07} & 2.47\e{-05} & 9.86\e{-13} & 3.19\e{-01}  \\
512 & 1.03\e{-02} & 5.33\e{-07} & 2.44\e{-05} & 1.50\e{-12} & 3.54\e{-01}  \\
1\sep,024 & 1.48\e{-02} & 5.81\e{-07} & 2.42\e{-05} & 2.34\e{-12} & 5.19\e{-01}  \\
2\sep,048 & 1.60\e{-02} & 5.90\e{-07} & 2.49\e{-05} & 5.48\e{-12} & 5.66\e{-01}  \\
4\sep,096 & 2.15\e{-02} & 5.70\e{-07} & 2.30\e{-05} & 1.39\e{-11} & 7.66\e{-01}  \\
8\sep,192 & 2.75\e{-02} & 5.69\e{-07} & 2.34\e{-05} & 1.71\e{-11} & 9.89\e{-01}  \\
16\sep,384 & 2.92\e{-02} & 5.66\e{-07} & 2.36\e{-05} & 2.71\e{-11} & 1.06\e{+00}   \\
32\sep,768 & 3.60\e{-02} & 1.03\e{-06} & 2.31\e{-05} & 9.68\e{-11} & 1.31\e{+00}   \\
65\sep,536 & 4.37\e{-02} & 5.36\e{-07} & 2.30\e{-05} & 2.31\e{-10} & 1.60\e{+00}   \\
131\sep,072 & 4.59\e{-02} & 5.65\e{-07} & 2.32\e{-05} & 4.64\e{-10} & 1.69\e{+00}   \\
262\sep,144 & 5.45\e{-02} & 5.66\e{-07} & 2.41\e{-05} & 6.96\e{-10} & 2.01\e{+00}   \\
524\sep,288 & 5.69\e{-02} & 6.03\e{-07} & 2.24\e{-05} & 1.58\e{-09} & 2.11\e{+00}   \\
1\sep,048\sep,576 & 6.67\e{-02} & 5.68\e{-07} & 2.42\e{-05} & 1.88\e{-09} & 2.46\e{+00}   \\
2\sep,097\sep,152 & 8.07\e{-02} & 5.86\e{-07} & 2.46\e{-05} & 6.20\e{-09} & 2.84\e{+00}   \\
4\sep,194\sep,304 & 7.91\e{-02} & 5.64\e{-07} & 2.40\e{-05} & 9.51\e{-09} & 2.97\e{+00}   \\
8\sep,388\sep,608 & 9.00\e{-02} & 8.71\e{-07} & 2.37\e{-05} & 1.76\e{-08} & 3.38\e{+00}   \\
16\sep,777\sep,216 & 1.01\e{-01} & 5.86\e{-07} & 2.61\e{-05} & 3.65\e{-08} & 3.81\e{+00}   \\
33\sep,554\sep,432 & 1.11\e{-01} & 5.80\e{-07} & 2.35\e{-05} & 7.77\e{-08} & 3.97\e{+00}   \\
67\sep,108\sep,864 & 1.17\e{-01} & 5.99\e{-07} & 2.36\e{-05} & 2.01\e{-07} & 4.43\e{+00}   \\
134\sep,217\sep,728 & 1.36\e{-01} & 6.25\e{-07} & 2.42\e{-05} & 3.74\e{-07} & 4.93\e{+00}   \\
 \bottomrule
 \end{tabular}

\end{center}
\vskip .5em
\caption{A comparison of the time required to evaluate Jacobi polynomials
via the algorithm of this paper with the time required to do so using
certain asymptotic expansions.  Here, the parameters in Jacobi's differential
equation were taken to be $a=-1/4$ and $b=1/3$.  All times are in seconds.}
\label{table:asym}
\end{table}

The asymptotic complexity of the procedure of 
Section~\ref{section:precomp} is
$\mathcal{O}\left(\log^2\left(N_{\mbox{\tiny max}}\right)\right)$;
however, the running times shown on the left-hand side of 
Figure~\ref{figure:asym} appear to grow more slowly than this.  This suggests that
a lower order term with a large constant is present.

We also compared the method of this paper with the results obtained
by using the well-known three-term recurrence relations satisfied
by solutions of Jacobi's differential equation
(which can be found in Section~10.8 of \cite{HTFII}, among many other sources)
 to evaluate the Jacobi polynomials.  Obviously,  such an approach is unsuitable as a mechanism
for evaluating a single Jacobi polynomial of large degree.  However, up to a certain point,
the recurrence relations are efficient and effective and it is of interest
to compare the approaches.  Table~\ref{table:rec} does so.  In these experiments,
$a$ was taken to be $1/4$  and $b$ was $-1/3$.

\begin{table}[h!]
\small
\begin{center}
 \begin{tabular}{lcccc}
 \toprule
 $N_{\mbox{\tiny max}}$ & Phase function  & Average evaluation time & Average evaluation time & Largest absolute \\
 & construction time & algorithm of this paper& recurrence relations       & error\\
 \midrule
32 & 2.63\e{-03} & 4.37\e{-07} & 1.74\e{-06} & 3.34\e{-13}  \\
64 & 2.56\e{-03} & 5.04\e{-07} & 2.33\e{-06} & 6.58\e{-13}  \\
128 & 5.65\e{-03} & 7.73\e{-07} & 3.72\e{-06} & 1.02\e{-12}  \\
256 & 9.84\e{-03} & 5.49\e{-07} & 6.16\e{-06} & 1.43\e{-12}  \\
512 & 1.02\e{-02} & 5.60\e{-07} & 1.15\e{-05} & 1.69\e{-12}  \\
1\sep,024 & 1.48\e{-02} & 5.82\e{-07} & 2.05\e{-05} & 8.31\e{-12}  \\
2\sep,048 & 1.59\e{-02} & 5.79\e{-07} & 3.76\e{-05} & 1.58\e{-11}  \\
4\sep,096 & 2.12\e{-02} & 5.63\e{-07} & 7.58\e{-05} & 3.83\e{-11}  \\
8\sep,192 & 2.73\e{-02} & 5.64\e{-07} & 1.51\e{-04} & 1.84\e{-10}  \\
16\sep,384 & 2.90\e{-02} & 5.68\e{-07} & 2.94\e{-04} & 1.98\e{-10}  \\
32\sep,768 & 3.58\e{-02} & 1.03\e{-06} & 5.96\e{-04} & 7.62\e{-11}  \\
 \bottomrule
 \end{tabular}

\end{center}
\vskip .5em
\caption{A comparison of the time required to evaluate Jacobi polynomials
via the algorithm of this paper with the time required to do so using
the well-known three-term recurrence relations.  Here, the parameters in Jacobi's differential
equation were taken to be $a=1/4$ and $b=-1/3$.  All times are in seconds.}
\label{table:rec}
\end{table}
\label{section:experiments:polys}

\end{subsection}


\begin{subsection}{Calculation of Gauss-Jacobi quadrature rules}
In this section, we describe experiments conducted to
measure the speed and accuracy with which
the algorithm of Section~\ref{section:quad} constructs Gauss-Jacobi
quadrature rules.  We used the  Hale-Townsend algorithm \cite{Hale-Townsend}, 
which appears to be the current state-of-the-art method for the numerical 
calculation of such rules,
as a basis for comparison.  It takes $\mathcal{O}\left(n\right)$ operations
to construct an $n$-point Gauss-Jacobi quadrature rule,
and calculates the quadrature nodes with
  double precision absolute accuracy and the quadrature weights with
double precision relative accuracy.  The Hale-Townsend
algorithm is faster and more accurate (albeit less general)
than the  Glaser-Liu-Rokhlin algorithm \cite{Glaser-Rokhlin}, which
can  be also used to construct $n$-point Gauss-Jacobi quadrature rules in $O(n)$ operations.
We note that the algorithm of
\cite{Bogaert} 
 for the construction  of Gauss-Legendre quadrature rules  (and not more general
Gauss-Jacobi quadrature rules)
is more efficient than the method of this paper.

\begin{figure}[h!!]
\begin{center}
\hfil
\includegraphics[width=.49\textwidth]{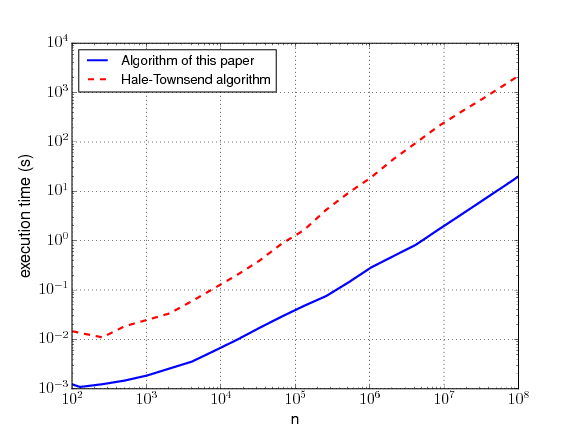}
\hfil 
\includegraphics[width=.49\textwidth]{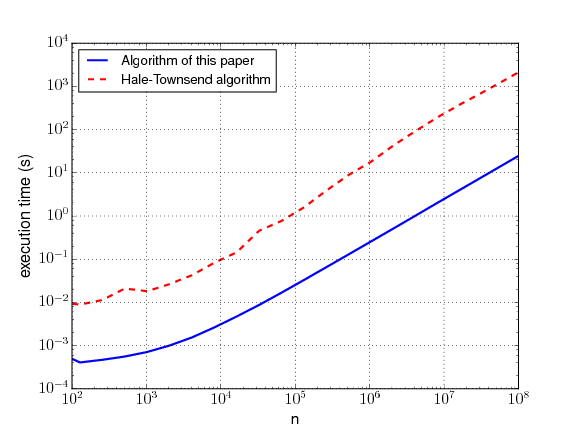}
\hfil
\end{center}
\caption{A comparison of the running time of the
algorithm of Section~\ref{section:quad} for the calculation
of Gauss-Jacobi quadrature rules with the algorithm of 
Hale-Townsend \cite{Hale-Townsend}.    In the experiments
whose results are shown on the left,
the parameters were taken to be 
 $a=0.25$ and $b=0.40$ while in the experiments corresponding
to the plot  on the right, the parameters were $a=-0.49$ and $b=0.25$.
} 
\label{figure:quad}
\end{figure}

For $a=0$ and $b=-4/10$ and various values of $n$,
we constructed the $n$-point Gauss-Jacobi quadrature rule 
using both the algorithm of Section~\ref{section:quad} and 
the Julia implementation \cite{Hale-Townsend-Code}
of the Hale-Townsend algorithm  made available
by the authors of \cite{Hale-Townsend}.
For each chosen value of $n$, we measured the relative difference in $100$ randomly selected
weights.  Unfortunately, in some cases
\cite{Hale-Townsend-Code} loses a small amount of accuracy when 
evaluating  quadrature weights  corresponding to quadrature nodes near the points $\pm 1$.
Accordingly,   when choosing random nodes,  we omitted those
corresponding to the $20$ quadrature nodes closest
to each of the points $\pm 1$.  We note that the loss
of accuracy in \cite{Hale-Townsend-Code}
 is merely a problem with that particular implementation of the Hale-Townsend
algorithm and not with the underlying scheme.  Indeed,
the MATLAB implementation of the Hale-Townsend algorithm 
included with the Chebfun package \cite{Chebfun}  does not suffer from this defect.  
We did not compare against the MATLAB version of the
code because it is somewhat slower than the Julia implementation.
Table~\ref{table:quad} reports the results.
We began our comparison with $n=101$ because 
when $n \leq 100$,
 the Hale-Townsend code combines
the well-known three-term recurrence relations satisfied by solutions
of Jacobi's differential equation
and Newton's method to construct the $n$-point Gauss-Jacobi quadrature rule.
This is also the strategy we recommend for constructing Gauss-Jacobi rules when
$n$ is small.

For different pairs of the parameters $a$ and $b$ and various values
of $n$, we used
 \cite{Hale-Townsend-Code} and the algorithm of this paper
to produce $n$-point Gauss-Jacobi quadrature rules and compared
the running times of these two algorithms.
Figure~\ref{figure:quad} displays the results.

\begin{table}[h!]
\small
\begin{center}
 \begin{tabular}{lcccc}
 \toprule
 $n$   & Running time of the & Running time of  & Ratio & Maximum relative\\
       & algorithm of Section~\ref{section:quad} & the algorithm of  \cite{Hale-Townsend} && error in weights  \\
 \midrule
101 & 4.45\e{-04} & 9.62\e{-03} & 2.15\e{+01}  & 4.47\e{-15}  \\
128 & 4.00\e{-04} & 9.26\e{-03} & 2.31\e{+01}  & 4.78\e{-15}  \\
256 & 4.54\e{-04} & 1.15\e{-02} & 2.54\e{+01}  & 5.76\e{-15}  \\
512 & 5.34\e{-04} & 2.14\e{-02} & 4.00\e{+01}  & 7.04\e{-15}  \\
1\sep,024 & 6.64\e{-04} & 2.27\e{-02} & 3.41\e{+01}  & 6.26\e{-15}  \\
2\sep,048 & 8.86\e{-04} & 3.00\e{-02} & 3.39\e{+01}  & 6.99\e{-15}  \\
4\sep,096 & 1.31\e{-03} & 5.28\e{-02} & 4.01\e{+01}  & 7.45\e{-15}  \\
8\sep,192 & 2.15\e{-03} & 8.76\e{-02} & 4.06\e{+01}  & 7.99\e{-15}  \\
16\sep,384 & 3.80\e{-03} & 2.68\e{-01} & 7.07\e{+01}  & 1.07\e{-14}  \\
32\sep,768 & 7.03\e{-03} & 4.45\e{-01} & 6.33\e{+01}  & 8.29\e{-15}  \\
65\sep,536 & 1.35\e{-02} & 7.91\e{-01} & 5.85\e{+01}  & 9.23\e{-15}  \\
131\sep,072 & 2.64\e{-02} & 1.61\e{+00}  & 6.10\e{+01}  & 1.04\e{-14}  \\
262\sep,144 & 5.22\e{-02} & 4.14\e{+00}  & 7.92\e{+01}  & 8.44\e{-15}  \\
524\sep,288 & 1.03\e{-01} & 9.06\e{+00}  & 8.74\e{+01}  & 1.09\e{-14}  \\
1\sep,048\sep,576 & 2.06\e{-01} & 1.80\e{+01}  & 8.73\e{+01}  & 1.29\e{-14}  \\
2\sep,097\sep,152 & 4.11\e{-01} & 4.26\e{+01}  & 1.03\e{+02}  & 1.26\e{-14}  \\
4\sep,194\sep,304 & 8.24\e{-01} & 9.27\e{+01}  & 1.12\e{+02}  & 1.16\e{-14}  \\
8\sep,388\sep,608 & 1.65\e{+00}  & 1.95\e{+02}  & 1.17\e{+02}  & 1.36\e{-14}  \\
16\sep,777\sep,216 & 3.30\e{+00}  & 3.86\e{+02}  & 1.16\e{+02}  & 1.43\e{-14}  \\
33\sep,554\sep,432 & 6.59\e{+00}  & 7.33\e{+02}  & 1.11\e{+02}  & 1.43\e{-14}  \\
67\sep,108\sep,864 & 1.31\e{+01}  & 1.42\e{+03}  & 1.08\e{+02}  & 1.48\e{-14}  \\
100\sep,000\sep,000 & 1.96\e{+01}  & 2.10\e{+03}  & 1.07\e{+02}  & 1.77\e{-14}  \\
 \bottomrule
 \end{tabular}

\end{center}
\vskip .5em
\caption{The results of an experiment comparing  the algorithm of \cite{Hale-Townsend}
for the numerical calculation of Gauss-Jacobi quadrature rules
with that of Section~\ref{section:quad}.  In these experiments,
 the parameters were taken to be $a=0$ and $b=-4/10$.  All times are in seconds.}
\label{table:quad}
\end{table}

\end{subsection}

\begin{subsection}{The Jacobi transform}
In these experiments, we measured the speed and accuracy of 
the algorithm of Section~\ref{section:transform} for the application
of the Jacobi transform and its inverse.

We did so in part by comparison with the  algorithm of Slevinsky \cite{Slevinsky1}
for applying the Chebyshev-Jacobi and Jacobi-Chebyshev transforms.
The Chebyshev-Jacobi transform
is the map which takes the coefficients of the Chebyshev expansion
of a function to the coefficients in its Jacobi expansion and the Jacobi-Chebyshev
transform is the inverse of this map.    
Although these are not the transforms we apply, the Jacobi transform
can be implemented easily by combining the method of
 \cite{Slevinsky1} with the nonuniform fast Fourier transform
(see, for instance, \cite{TownsendLegendre}  and \cite{TownsendLegendre2} 
for an approach of this type for applying the Legendre transform and its inverse).
Other methods for applying the Jacobi transform, some of which have
 lower asymptotic complexity than the algorithm
of \cite{Slevinsky1} and the approach of this paper, are available.
Butterfly algorithms  such as \cite{Li-Yang-Martin-Ho-Ying,Li-Yang-Ying1,Oneil-Rokhlin}
allow for the application of the Jacobi transform and various related mappings
 in $\mathcal{O}\left(n\log(n)\right)$ operations; however, existing
methods  are either numerically unstable or they require 
expensive precomputation phases with higher asymptotic complexity.
The Alpert-Rokhlin method \cite{Alpert-Rokhlin} uses a multipole-like
approach to apply the Chebyshev-Legendre and Legendre-Chebyshev transforms
in $\mathcal{O}\left(n\right)$ operations.  The Legendre transform can be
computed in $\mathcal{O}\left(n\log(n)\right)$
operations by combining
this algorithm with the fast Fourier transform.
This approach is extended in  \cite{Keiner} in order to compute
expansions in Gegenbauer polynomials in $\mathcal{O}(n\log(n))$ operations.
It seems likely that these methods can be extended to compute
expansions in more general classes of Jacobi polynomials; however,
to the authors' knowledge no such algorithm has been published 
and algorithms from this class require expensive precomputations.
A further discussion of methods for the application of the Jacobi transform
and related mappings can be found in \cite{Slevinsky1}.

\begin{table}[t!]
\small
\begin{center}
 \begin{tabular}{lccc}
 \toprule
 $n$   & Forward Jacobi   & Chebyshev-Jacobi             & Ratio  \\
       & transform time    & transform time &  \\
       & algorithm of Section~\ref{section:transform} & algorithm of \cite{Slevinsky1} &   \\
 \midrule
10 & 6.78\e{-07} & 1.96\e{-04} & 2.88\e{+02}   \\
16 & 5.10\e{-07} & 1.78\e{-04} & 3.49\e{+02}   \\
32 & 2.21\e{-05} & 1.99\e{-04} & 9.00\e{+00}   \\
64 & 3.08\e{-05} & 3.86\e{-04} & 1.25\e{+01}   \\
128 & 2.74\e{-05} & 5.65\e{-04} & 2.05\e{+01}   \\
256 & 4.52\e{-05} & 8.33\e{-04} & 1.84\e{+01}   \\
512 & 8.46\e{-05} & 2.77\e{-03} & 3.27\e{+01}   \\
1\sep,024 & 1.56\e{-04} & 5.71\e{-03} & 3.64\e{+01}   \\
2\sep,048 & 3.44\e{-04} & 1.53\e{-02} & 4.46\e{+01}   \\
4\sep,096 & 9.41\e{-04} & 2.75\e{-02} & 2.92\e{+01}   \\
8\sep,192 & 2.35\e{-03} & 8.75\e{-02} & 3.70\e{+01}   \\
16\sep,384 & 6.08\e{-03} & 2.12\e{-01} & 3.48\e{+01}   \\
32\sep,768 & 1.43\e{-02} & 5.08\e{-01} & 3.55\e{+01}   \\
65\sep,536 & 2.92\e{-02} & 9.91\e{-01} & 3.39\e{+01}   \\
131\sep,072 & 7.36\e{-02} & 3.25\e{+00}  & 4.42\e{+01}   \\
262\sep,144 & 1.35\e{-01} & 4.47\e{+00}  & 3.30\e{+01}   \\
524\sep,288 & 4.50\e{-01} & 1.54\e{+01}  & 3.43\e{+01}   \\
1\sep,048\sep,576 & 1.05\e{+00}  & 2.28\e{+01}  & 2.15\e{+01}   \\
2\sep,097\sep,152 & 2.72\e{+00}  & 6.70\e{+01}  & 2.45\e{+01}   \\
4\sep,194\sep,304 & 8.55\e{+00}  & 1.44\e{+02}  & 1.69\e{+01}   \\
8\sep,388\sep,608 & 1.74\e{+01}  & 3.33\e{+02}  & 1.91\e{+01}   \\
16\sep,777\sep,216 & 3.67\e{+01}  & 5.80\e{+02}  & 1.57\e{+01}   \\
33\sep,554\sep,432 & 7.92\e{+01}  & 1.61\e{+03}  & 2.03\e{+01}   \\
67\sep,108\sep,864 & 1.66\e{+02}  & 3.03\e{+03}  & 1.81\e{+01}   \\
100\sep,000\sep,000 & 1.52\e{+02}  & 5.25\e{+03}  & 3.43\e{+01}   \\
 \bottomrule
 \end{tabular}

\end{center}
\vskip .5em
\caption{
A comparison of the time required to apply the
Chebyshev-Jacobi transform  using the algorithm
of \cite{Slevinsky1} with the time required to apply the forward Jacobi
transform via the algorithm of Section~\ref{section:transform}.
Here, $a$ was taken to be $1/4$ and $b$ was $-4/10$.  All times are in seconds. 
}
\label{table:transform}
\end{table}

\begin{figure}[h!!]
\begin{center}
\hfil
\includegraphics[width=.49\textwidth]{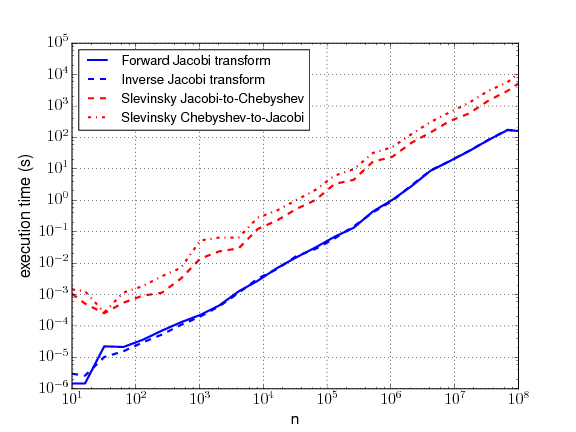}
\hfil 
\includegraphics[width=.49\textwidth]{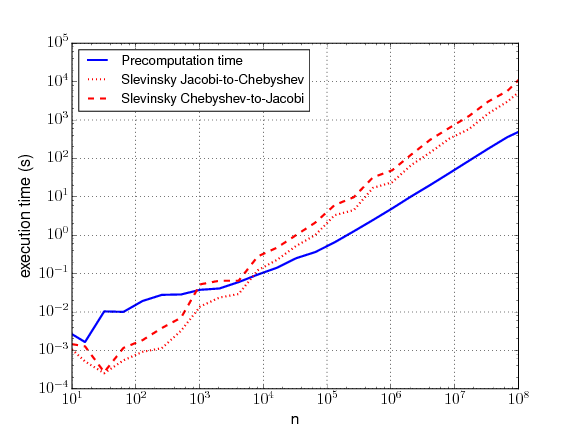}
\hfil
\end{center}
\caption{On the left is a comparison of the time required to apply the
Chebyshev-Jacobi transform and its inverse using the algorithm
of \cite{Slevinsky1} with the time required to apply the forward and inverse
Jacobi transforms via the algorithm of Section~\ref{section:transform}.
On the right is a comparison of the cost of the precomputations
for the procedure of Section~\ref{section:transform} with the time required
to apply the Jacobi-Chebyshev map and its inverse via the 
algorithm of \cite{Slevinsky1}.
In these experiments, $a$ was taken to be $-1/4$ and $b$ was $0$.
}
\label{figure:transform1}
\end{figure}

Figure~\ref{figure:transform1} and Table~\ref{table:transform}
report the results of our comparisons
with the Julia implementation \cite{Slevinsky-code}
of Slevinsky's algorithm.  The graph on the left side of Figure~\ref{figure:transform1}
compares the time taken by the algorithm of Section~\ref{section:transform}
to apply the Jacobi transform and its inverse
with the time required to apply the Chebyshev-Jacobi mapping and its
inverse via \cite{Slevinsky1}, while the graph on the right 
 compares the time required by our precomputation procedure 
with the time required to apply the 
Chebyshev-Jacobi mapping and it inverse with Slevinsky's algorithm.
We observe that cost of our algorithm, including the precomputation stage, is less
than that of \cite{Slevinsky1} at relatively modest orders.
Moreover, Figure~\ref{figure:transform1} strongly suggests that
the asymptotic running time of our algorithm for the application
of the Jacobi transform is similar to the
$O\left(n \log^2\left(n\right) / \log\left(\log\left(n\right)\right)\right)$
complexity of Slevinsky's algorithm.

Owing to the loss of accuracy which arises when the Formula~(\ref{introduction:ptilde}) 
is used to evaluate Jacobi polynomials of large degrees,
we expect the error in the Jacobi transform
of Section~\ref{section:transform} to increase
as the order of the transform increases.  This is indeed the case,
at least when it is applied to functions whose Jacobi coefficients
do not decay or decay slowly.  However, when the transform is applied
to   smooth functions, whose Jacobi coefficients decay rapidly,
the errors grow more slowly than in the general case.
  This is the same as the behavior 
of the algorithms \cite{Slevinsky1} and \cite{TownsendLegendre}.
We carried out a further set of experiments to illustrate this effect.
In particular,  we applied the forward Jacobi transform followed
by the inverse Jacobi transform to vectors which decay at various rates.
We constructed test vectors by choosing their entries
from a Gaussian distribution and then scaling them so as to
achieve a desired rate of decay.
Figure~\ref{figure:transform2} reports the results.
It also contains a plot of the rank of the matrix (\ref{transform:matrix})
as a function of $n$ for various pairs of the parameters $a$ and $b$.

\begin{figure}[t!!]
\begin{center}
\hfil
\includegraphics[width=.49\textwidth]{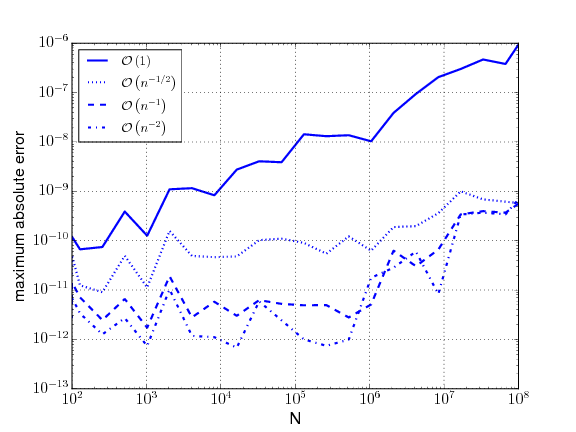}
\hfil
\includegraphics[width=.49\textwidth]{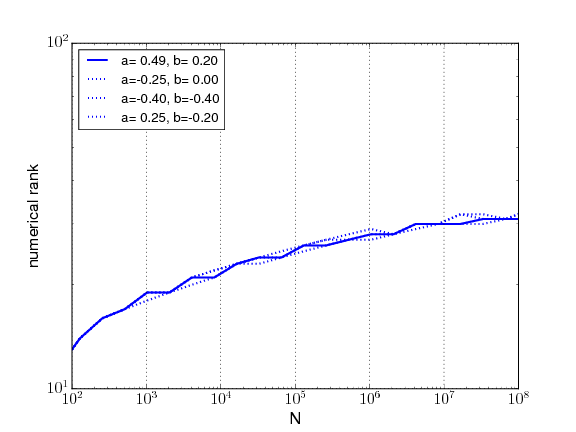}
\hfil
\end{center}
\caption{On the left are plots of the largest absolute errors which occurs 
when the composition of the inverse
and forward Jacobi transforms of Section~\ref{section:transform} are applied to vectors whose coefficients
decay at varying rates.   For these experiments, $a=-1/4$ and $b=0$.
 On the right are plots of the rank of the matrix
(\ref{transform:matrix}) 
as a function of $N_{\mbox{\tiny max}}$ for various values of the parameters $a$ and $b$.
}
\label{figure:transform2}
\end{figure}

We also compared the time required to apply the forward Jacobi transform
via the algorithm of Section~\ref{section:transform} with
the time required to do so by evaluating the matrix $\mathcal{J}_n^{(a,b)}$
using the well-known three-term recurrence relations and then
applying it directly (we refer to this as the
``brute force technique'').  This is the methodology we recommend for
transforms of small orders. In these experiments,
the parameters $a$ and $b$ were taken to be $a=1/4$ and $b=-1/3$.
  Figure~\ref{figure:transform3}  shows the results. 

\begin{figure}[t!!]
\begin{center}
\hfil
\includegraphics[width=.49\textwidth]{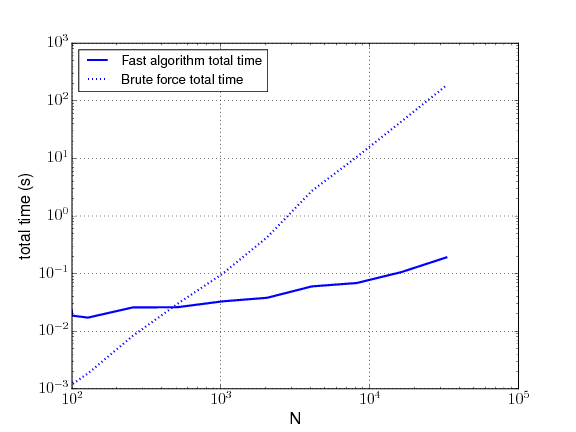}
\hfil
\includegraphics[width=.49\textwidth]{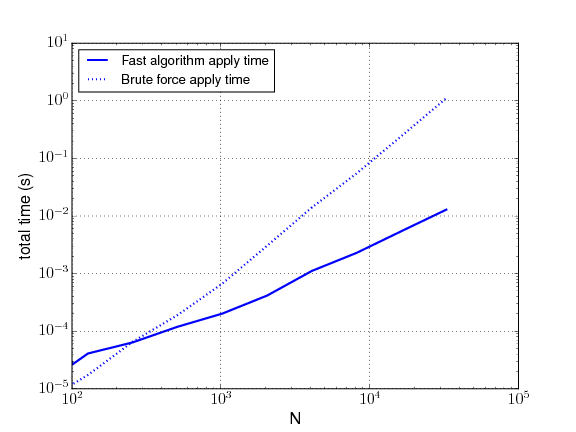}
\hfil
\end{center}
\caption{A comparison of the time take to apply the Forward Jacobi
transform via ``brute force'' with the time required to do so
via the algorithm of Section~\ref{section:transform}.
On the left is a plot of the total time taken.  For the algorithm
of Section~\ref{section:transform} this includes the time taken
by the precomputation phase, while for ``brute force'' technique this includes the time required
to evaluate the entries of the matrix $\mathcal{J}_n^{(a,b)}$
via three-term recurrence relations.
On the right is a comparison of the application times only.
Here, the parameters are $a=1/4$ and $b=-1/3$.
}
\label{figure:transform3}
\end{figure}


\label{section:experiments:transform}
\end{subsection}

\label{section:experiments}
\end{section}

\begin{section}{Conclusion and Further Work}

We have described a suite of fast algorithms for forming and manipulating expansions
in Jacobi polynomials.     They are based on the well-known fact that
Jacobi's differential equation admits a nonoscillatory phase function.
Our algorithms use  numerical methods, rather than asymptotic expansions, to evaluate
the phase function and other functions related to it.  
We do so in part because existing asymptotic expansions for the phase 
function do not suffice for our purposes (they either not sufficiently
accurate or they are not numerically viable), but also because such
techniques can be easily applied to any family of special functions
satisfying a second order differential equation with nonoscillatory
coefficients.  We will report on the 
application of our methods to other families of special functions
at a later date.

It would be of some interest to accelerate the procedure of Section~\ref{section:precomp}
for the construction of the nonoscillatory phase and amplitude functions.
There are a number obvious mechanisms for doing so, but perhaps the most
promising is the observation that the ranks of matrices with entries
\begin{equation}
\psi^{(a,b)}(t_j,\nu_k)
\end{equation}
and
\begin{equation}
M^{(a,b)}(t_j,\nu_k)
\end{equation}
are quite low --- indeed, in the experiments of  this paper
they were never observed to be larger than $40$.
 This means that, at least in principle, the nonoscillatory
phase and amplitude can be represented
via $40 \times 40$ matrices, and that a carefully designed spectral scheme
which takes this into account 
could compute the $40$ required 
solutions of the  ordinary differential equation (\ref{phase:Neq})
extremely efficiently.

The authors are investigating such an approach to the construction
of $\psi^{(a,b)}$ and $M^{(a,b)}$.

As the parameters $a$ and $b$ increase beyond $\frac{1}{2}$, our algorithms become less
accurate, and they ultimately fail.  This happens because
for values of the parameters $a$ and $b$ greater than $\frac{1}{2}$,
 the  Jacobi polynomials have turning points and 
the crude approximation (\ref{introduction:psiapprox}) becomes inadequate.
An obvious remedy is to use a more sophisticated approximations
for $\psi^{(a,b)}$ and $M^{(a,b)}$.  
The authors will report on extensions of this work which make use
of such an approach at a later date.

\label{section:conclusion}
\end{section}

\begin{section}{References}
\bibliographystyle{acm}
\bibliography{jacobi}
\end{section}

\end{document}